# FINDING EIGENVECTORS: FAST AND NONTRADITIONAL APPROACH

UDITA N. KATUGAMPOLA

Abstract. Diagonalizing a matrix $A$, that is finding two matrices $P$ and $D$ such that $A = PDP^{-1}$ with $D$ being a diagonal matrix needs two steps: first find the eigenvalues and then find the corresponding eigenvectors. We show that we do not need the second step when diagonalizing matrices with a spectrum, $|\sigma(A)| \leq 2$ since those vectors already appear as nonzero columns of the eigenmatrices, a term defined in this work. We further generalize this for matrices with $|\sigma(A)| > 2$ and show that eigenvectors lie in the column spaces of eigenmatrices of the complementary eigenvalues, an approach without using the classical Gaussian elimination of rows of a matrix. We also provide several shortcut formulas to find eigenvectors that does not use echelon forms. We may summarizes this new method with the mnemonic "Find your puppy at your neighbors'!".

*Dedicated to*

අම්මා *and* තාත්තා

## 1. Introduction

For decades if not centuries, we have been overlooking the use of columns when it comes to producing eigenvectors, even though they are sitting in the column space of a matrix. The standard approach is to find a basis for the null$(\mathbf{A} - \lambda\mathbf{I})$ and take the basis elements as the eigenvectors associated with $\lambda$. One way to do that is to use the Gaussian Elimination process to row-reduce the matrix $\mathbf{A} - \lambda\mathbf{I}$ and solve for the nonzero vector(s), $\nu$ that makes $(\mathbf{A} - \lambda\mathbf{I})\nu = 0$. However, it is an easy exercise to show that an eigenvector is indeed a linear combination of columns when 0 is not the associated eigenvalue. If an eigenvalue $\lambda \neq 0$, we have by the definition of an eigenvector

$$\mathbf{A}\nu = \lambda\nu \tag{1.1}$$

Taking $\mathbf{A} = [a_1\ a_2\ \cdots\ a_n]$ and $\nu = [\nu_1,\ \nu_2,\ \cdots,\ \nu_n]^T$, we may rewrite (1.1) in the form

$$\nu = \frac{\nu_1}{\lambda}a_1 + \frac{\nu_2}{\lambda}a_2 + \cdots + \frac{\nu_3}{\lambda}a_3$$

Thus, the eigenvector itself is a linear combination of columns of $\mathbf{A}$. On the other hand, if 0 is an eigenvalue, the corresponding eigenvector is in the null $\mathbf{A}$ and there is apparently no connection to the column space of $A$, other than the components of the eigenvector, in fact, being the weights that makes the corresponding linear combination 0, sending $\nu$ to the nullspace of $A$.

Now, suppose there is another eigenvalue $\mu \neq 0$ to $\mathbf{A}$. Then, it can be shown that the eigenvector $\nu$ corresponding to 0 lies in $\mathrm{col}\,(\mathbf{A} - \mu\mathbf{I})$, instead. For example, let $\mathbf{A} = \begin{pmatrix} -1 & 1 \\ -2 & 2 \end{pmatrix}$, which has eigenvalues $\lambda = 0, 1$ with associated eigenvectors $v_0 = [1\ \ 1]^T$ and $v_1 = [1\ \ 2]^T$, respectively. It is easy to see that the eigenvector of 0, $v_0 \in \mathrm{col}(\mathbf{A} - 1 \cdot \mathbf{I})$. In addition, the eigenvector of 1, $v_1 \in \mathrm{col}(\mathbf{A} - 0 \cdot \mathbf{I}) = \mathrm{col}(\mathbf{A})$. In other words, the nonzero columns of certain matrices or their linear combinations are the eigenvectors. The goal is to find that linear combination.

Before diving further into the details, let's have a look at some historical development of the subject. We have adopted several historical remarks from five major resources that includes [3, 32, 34, 46, 50] and sometimes have taken verbatim to keep the ideas of the original authors intact.

Among many others who contributed to the development of spectral theory, several names are note worthy. This does not imply in anyway that the other works are not important, but proceed in this Department of Applied Mathematics, Florida Polytechnic University, Lakeland FL 33805, USA..

Email: uditanalin@yahoo.com

manner due to space-time limitations. Interested readers are referred to the sources in question or others for further historical information.

The groundbreaking work laid by Daniel Bernoulli (1700-1784) is considered to be the first step in the direction of spectral theory as it is known today, even though there were no explicit mention of eigenvalues, eigenvectors or any similar terminology at that time. However, he used the idea behind the eigenvectors (eigenfunctions) to represent the vibrating string (rope) as the limit of $n$ oscillating particles and considered eigenvectors as configurations in his work in "*Theoremata de oscillationibus corporum filo flexili connaxorum et catenae verticaliter suspensae*"(1732) [6] as illustrated in Figure 1. Even though there was no mention about the naming convention the idea of eigenvalues (Figure 2) appeared in Cauchy's work in "*Sur l'équation à l'aide de laquelle on détermine les inégalités séculaires des mouvements des planétes*"(1829) as shown in Figure 2 and he referred to the eigenvalues as *valeurs propres* (proper values). It also shows his unintentional use of the idea of matrices in writing the system of equation as an array without the brackets or parenthesis as we use today. and the name of *characteristic equation* explicitly appeared in "*Mémoire sur L'intégration des équations linéaires*"(1840) shown in Figure 3. The systematic introduction of the concept of matrices was a result of James J. Sylvester (1814-1897) and Arthur Cayley (1821-1895) in an attempt to find an algebraic description of $n$-dimensional spaces. In 1852, Sylvester showed that 'eigenvalues', as we know today, are the roots of the characteristic polynomial $\det(\lambda I - A)$ and also introduced the name *determinant* to the literature of mathematics. In 1858, Cayley inaugurated the theory of matrices, possibly motivated by Cauchy's use of rectangular array in Equation 11 in Figure 2, in which he also discussed a diagonalization process of matrices [46].

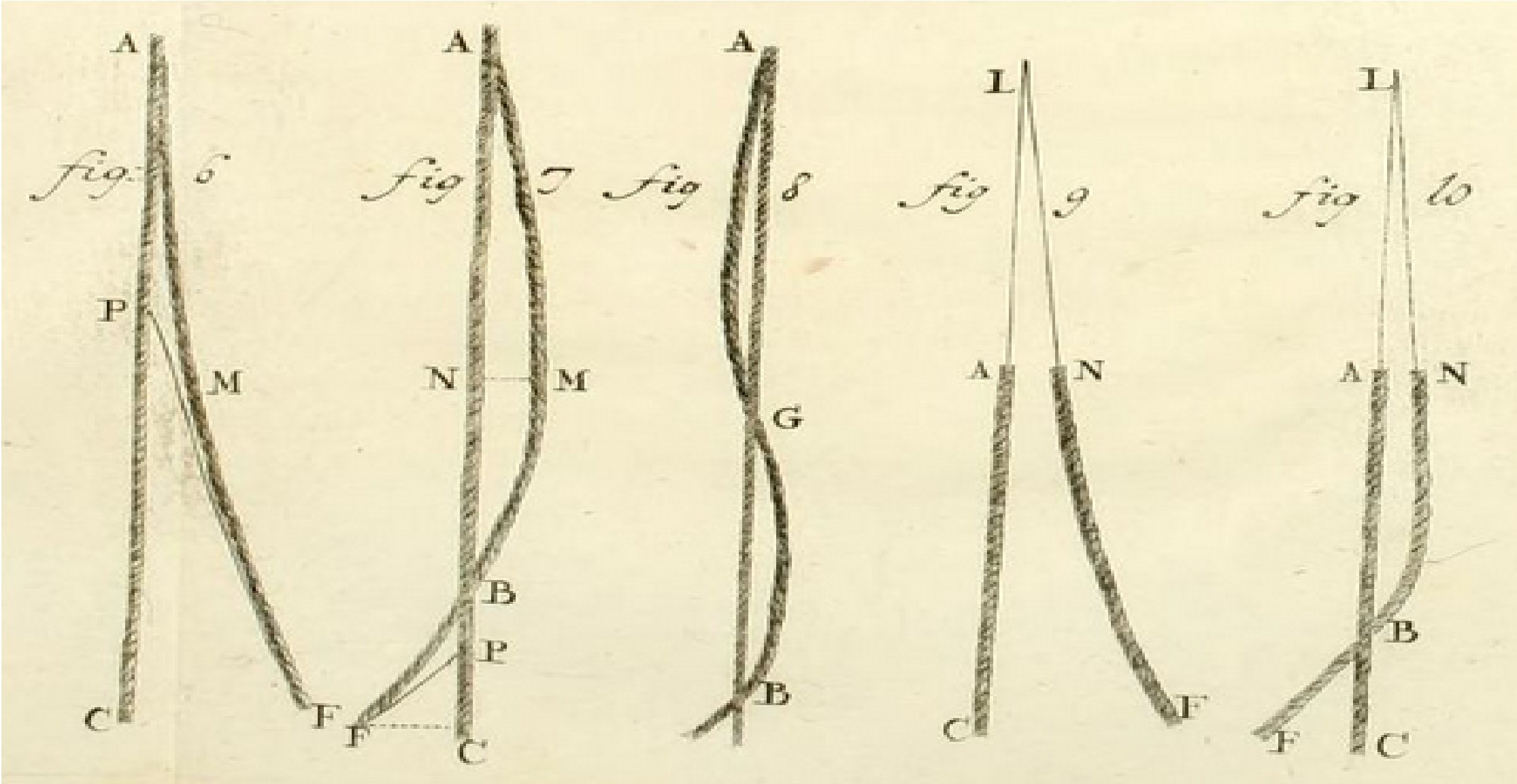

FIGURE 1. The three eigenvectors for a uniform 3-body system: Bernoulli's 1732 illustration [6, p. 108], from a scan of the journal at Archive.org

Some historical remarks about the uses of half-German names *eigenvalues* and *eigenvectors* are also in order. D. Bernoulli did not use any of those words in his work on vibrating strings in 1732 [6], even though, ironically he had used the letters $\lambda$ and $\mu$ that we use to represent eigenvalues today. As mentioned earlier, the notion of eigenvectors were not identified or discussed, even though, he used different configurations to represent a similar idea [6, Figure on p.108]. A similar argument appeared in Lecture notes of [17].

The eigenvalues pay a vital role in modern mathematics. However, the use of the name, a half-translation of the German *eigenwert*, was not settled until recently. In older works one may find latent root, characteristic value, secular (value) and proper value, as alternatives [27, 33].

According to the excellent expository and encyclopedic article about the word *eigen-* and its origin in mathematics maintained by J. Miller in [34], the term *characteristic* derived from Cauchy's 1840 work et pourront s'écrire comme il suit :

$$(10)\quad \begin{cases} (A_{xx}-s)x + A_{xy}y + A_{xz}z + \ldots = 0, \\ A_{xy}x + (A_{yy}-s)y + A_{yz}z + \ldots = 0, \\ A_{xz}x + A_{yz}y + (A_{zz}-s)z + \ldots = 0, \\ \ldots\ldots\ldots\ldots\ldots\ldots\ldots\ldots \end{cases}$$

Cela posé, il résulte des principes établis dans le Chapitre III de l'*Analyse algébrique* (§ 2) (¹) que le premier membre de l'équation (8), ou S, sera une fonction alternée des quantités comprises dans le Tableau :

$$(11)\quad \begin{cases} A_{xx}-s, & A_{xy}, & A_{xz}, & \ldots, \\ A_{xy}, & A_{yy}-s, & A_{yz}, & \ldots, \\ A_{xz}, & A_{yz}, & A_{zz}-s, & \ldots, \\ \ldots, & \ldots, & \ldots\ldots, & \ldots, \end{cases}$$

savoir celle dont les différents termes sont représentés, aux signes près, par les produits qu'on obtient, lorsqu'on multiplie ces quantités, $n$ à $n$, de toutes les manières possibles, en ayant soin de faire entrer dans chaque produit un facteur pris dans chacune des lignes horizontales du Tableau et un facteur pris dans chacune des lignes verticales. En opérant ainsi, on trouvera, par exemple, pour $n = 2$,

$$(12)\quad S = (A_{xx}-s)(A_{yy}-s) - A_{xy}^2;$$

Figure 2. Cauchy's use of eigenvalues in his work on planetary motion. Image courtesy of the Bibliothèque Nationale de France's Gallica collection [19]

D'ailleurs comme l'élimination des facteurs A, B, C, ... entre les formules (4), fournira une *équation caractéristique*

$$(5)\qquad \mathcal{S} = 0$$

qui sera du degré $n$ par rapport à $s$, la valeur de $\mathcal{S}$ étant

$$(6)\qquad \mathcal{S} = (s + \mathcal{L})(s + \mathcal{Q})\ldots - \mathfrak{M}\mathcal{P}\ldots\ldots + \text{etc.},$$

on pourra, dans les formules (3), prendre pour $s$ une quelconque des $n$ racines de l'équation (5). Il y a plus : comme, étant donnés, pour les variables principales, deux ou plusieurs systèmes de valeurs propres à vérifier les équations (1), on obtiendra de nouvelles intégrales de ces mêmes équations en ajoutant l'une à l'autre les diverses valeurs de chaque variable principale, il est clair qu'on vérifiera encore les équations (1) en in which he introduced the term *l'equation caractéristique* and investigated its roots in his *Mémoire sur l'integration des équations linéaires* [11]. Carl G. J. Jacobi used the word *secular* to mean the Cauchy's term, characteristic, in his famous paper *Über ein leichtes Verfahren die in der Theorie der Säcularstörungen vorkommenden Gleichungen numerisch aufzulösen* (1846) [28] that studied celestial mechanics and introduced his eigenvalue problem. The name then transformed to the *latent* by James J. Sylvester (1814-1897) in the 1883 paper *On the Equation to the Secular Inequalities in the Planetary Theory* Phil. Mag. **16**, 267.

Figure 3. Cauchy's explicit use of the name *characteristic equation.* Image courtesy of the Archive.org [4]

The *eigen* terminology is credited to David Hilbert (1862-1943), though he may have had some influence from contemporary mathematicians. *Eigenfunktion* and *Eigenwert* appeared in Hilbert's work on integral equations in *Grundzüge einer allgemeinen Theorie der linearen Integralgleichungen* (1904) [24]. Hilbert's original starting point was a non-homogeneous integral equation with a parameter $\lambda$ for which the matrix counterpart is $(I - \lambda A)x = y$. Hilbert called the values of $\lambda$ that generate nonzero solutions to the homogeneous version of these equations Eigenwerte; they are the reciprocals of the characteristic/latent roots of $A$, adding some confusion to the present meaning [34].

*Eigenwert* was used in a different context in J. von Neumann's work saying that *"Ein Eigenwert ist eine Zahl, zu der es eine Funktion $f \neq 0$ mit $Rf = \lambda f$ ist dann Eigenfunktion."* [37]. Hilbert's *eigen*-terminology was continually promoted and became the dominant usage afterward, as can be seen even today [34]. Finally the name *eigenvalue* became a synonymous for its predecessors such as characteristic value and latent root.

Other than this terminology another term rendered from *eigen*- has been used in some $19^{th}$ century English writings. The terms *proper values* and *proper function* appear in von Neumann's work, but defeated by the *eigen*- equivalents later in the literature [34].

As argued in [34], the terminology *Eigenvector* appears in R. Brauer & H. Weyl's *Spinors in n Dimensions*, Amer. J. Math. (1935). While *eigenvalue*, *eigenfunction* and *eigenvector* were translations of German originals, the eigenstate of Dirac's Principles of Quantum Mechanics (1930) [14, p.35] was a new addition to the philosophy and it marks the arrival of eigen as an English brainchild. *Eigenstate* was translated into German as *Eigenzustand* [34].

In this work we introduced a new terminology, the *eigenmatrix* to the mathematics literature as the lost ancestor/child of the sequence: eigenvalue and eigenvector, in that order, with the explanation that eigenmatrix produces eigenvalues and eigenvalues and eigenmatrices produce eigenvectors. We use $\kappa$ to represent eigenmatrices, $\lambda$ or $\mu$ for eigenvalues and $\nu$ for eigenvectors, whenever possible. Different from the row-reduction approach, there are other ways to find eigenvectors given the corresponding eigenvalues. A survey of eigenvalue-eigenvector duality for Hermitian matrices is presented in the recent work of [15].

## 2. Observations and Preparations

This work sheds some light on two ideas in spectral theory: finding the eigenvectors and diagonalizing matrices. We use the standard terminology in linear algebra, except in few cases, that can be found in any linear algebra book, for example, in Friedberg's [18] or Lay's [31]. We denote the spectrum of a matrix $A$, that is the set of eigenvalues, by $\sigma(A)$ and the matrices are taken over commutative rings such as $\mathbb{R}$ or $\mathbb{C}$, unless otherwise specified, and do not distinguish most of the time. Several fundamental results, without proofs, are in order to make the writing self-contained. Readers may skip the current section and jump on reading to the next section without hindering the continuity of the arguments present in this work.

The familiar approach to diagonalizing an $n \times n$ matrix $A$ needs two steps: first find the eigenvalues $\lambda_1, \dots, \lambda_n$ and then find the corresponding eigenvectors $v_1, \dots v_n$. If the $n$ eigenvectors are linearly independent, we form

$$P = (v_1 \ \dots \ v_n), \quad \text{and} \quad D = \begin{pmatrix} \lambda_1 & \cdots & 0 \\ \vdots & \ddots & \vdots \\ 0 & \dots & \lambda_n \end{pmatrix}, \tag{2.1}$$

this leads to

$$A = PDP^{-1}. \tag{2.2}$$

Eigenvectors appear as basis elements of the null-space (or the eigenspace of $\lambda_i$) of $A - \lambda_i I$ for each eigenvalue $\lambda_i$. In the classical approach we use row-operations to reduce the matrix $A - \lambda_i I$ to echelon form and find non-zero solutions to $A - \lambda_i I = \mathbf{0}$ as eigenvectors of $\lambda_i$. Before state the next result, for completeness, let's give the following definitions.

**Definition 2.1** (Eigenvalue and Eigenvector)**.** Let $A \in \mathbb{C}^{n^2}$ be an $n \times n$ matrix. The roots of the characteristic equation $\det(A - \lambda I) = 0$ are called the eigenvalues of $A$. If $\nu \neq 0$ is such that $(A - \lambda I)\,\nu = 0$, we call $\nu$ is an eigenvector of $A$ associated with the eigenvalue $\lambda$.

The set $E_\lambda = \{\nu \in \mathbb{C}^n \big| (A - \lambda I)\,\nu = 0\}$ is called the eigenspace of the eigenvalue $\lambda$. The number of basis elements of the eigenspace $E_\lambda$ is called the *geometric multiplicity* of $\lambda$.

Sometimes it is important to discuss fast approaches, or shortcut methods, when these concepts are introduced in the classroom mainly to use them in other applications. There is a common trick that we use for $2 \times 2$ matrices [30, 44, 52, 56].

Let $A = \begin{bmatrix} a & b \\ c & d \end{bmatrix}$ be a real matrix. We first find eigenvalues solving $\det(A - \lambda I) = 0$ for $\lambda$. This leads to the characteristic equation

$$\det(A - \lambda I) = (a - \lambda)(d - \lambda) - bc = \lambda^2 - \mathrm{Tr}(A)\lambda + \det(A) = 0, \tag{2.3}$$

where $\mathrm{Tr}(A) = a + d = \lambda_1 + \lambda_2$ is the *trace* of $A$ and $\det(A) = \lambda_1\lambda_2$, where $\lambda_1$ and $\lambda_2$ are the two eigenvalues of $A$ that can be found solving the quadratic equation (2.3) for $\lambda$. We do not distinguish between real or complex eigenvalues and the same argument works for both cases.

We then proceed as follows to find the corresponding eigenvectors. An eigenvector associated with the eigenvalue $\lambda_1$ is a nontrivial solution $\nu_1$ such that

$$(A - \lambda_1 I)\,\nu_1 = 0. \tag{2.4}$$

This leads to

$$A - \lambda_1 I = \begin{pmatrix} a - \lambda_1 & b \\ c & d - \lambda_1 \end{pmatrix} \tag{2.5}$$

The matrix $A - \lambda_1 I$ must be singular and should have rank 1. This, in turn, implies that the second row is a multiples of the first row. That means the Gaussian elimination for echelon form would leads to a row of zeros. Thus, all we need to find the eigenvector is the first row. In particular, if $\nu_1 = [v_1, v_2]^T$, then (2.4) implies

$$(a - \lambda_1)v_1 + bv_2 = 0.$$

To find a nonzero eigenvector for $\lambda_1$, we only need to find a particular set of numbers for $v_1$ and $v_2$ that would produce a nonzero $\nu_1$. That means at least one of $v_1$ and $v_2$ must be nonzero. In the case of $a - \lambda_1 = 0$ and $b \neq 0$, we may take $v_1 = 1$ and $v_2 = 0$. The case of $b = 0$ and $a - \lambda_1 \neq 0$ is similar. In the general case we may take, for example, $v_1 = -b$ and $v_2 = a - \lambda_1$. This would lead to the eigenvector

$$\nu_1 = \begin{bmatrix} -b \\ a - \lambda_1 \end{bmatrix}.$$

If both $a - \lambda_1$ and $b$ are zero, we may use the second row to find an eigenvector given by

$$\nu_1 = \begin{bmatrix} d - \lambda_1 \\ -c \end{bmatrix}.$$

The minus sign (-) can appear with any component of the vector since any nonzero scalar multiple of an eigenvector is still an eigenvector. Let's consider two examples to illustrate the result.

**Example 2.2.** Let

$$A = \begin{bmatrix} 3 & 1 \\ 2 & 4 \end{bmatrix}.$$

The characteristic polynomial is

$$\lambda^2 - (3 + 4)\lambda + ((3)(4) - (1)(2)) = \lambda^2 - 7\lambda + 10,$$

which leads to

$$\lambda = \frac{7 \pm \sqrt{49 - 4(10)}}{2} = 2, 5.$$

Let $\lambda_1 = 2$ and $\lambda_2 = 5$. Now for an eigenvector of $\lambda_1 = 2$

$$A - \lambda_1 I = \begin{pmatrix} 1 & 1 \\ 2 & 2 \end{pmatrix}$$

It is easy to see that the second row of the matrix in question is a multiple of the first, proving the claim. Now using the shortcut method above, we can find one eigenvector of $\lambda_1 = 2$ to be

$$\nu_1 = \begin{bmatrix} -1 \\ 1 \end{bmatrix}$$

A similar argument finds an eigenvector associated with the eigenvalue $\lambda_2 = 5$ given by $\nu_2 = [-1 \quad -2]^T$, which may be scaled and negate to produce $\nu_2 = [1, \quad 2]^T$. The method equally applies to matrices with complex eigenvalues as well. They will appear as conjugate pairs when the matrix is in $\mathbb{R}^{2\times 2}$.

We recall several results that are fundamental in Algebra and Linear Algebra and are heavily used throughout this work. These results can be found in any standard books in those areas. The first is the following, the proof of which can be found in any linear algebra text.

**Theorem 2.3** (Eigenspace). *Eigenspace is an invariant subspace.*

This provides a valuable tool in our work, specially when expressing an eigenvector as a linear combination of column vectors of eigenmatrices, a term that will be defined later in this section. Next we have one of the fundamental results in algebra, the so-called Fundamental Theorem of Algebra [16], which allows us to factor the characteristic polynomial into linear factors, the main tool in finding the eigenvectors as proposed in this work.

**Theorem 2.4** (Fundamental Theorem of Algebra). *The field $\mathbb{C}$ is algebraically closed.*

The result tells us that every non-constant single-variable polynomial with complex coefficients has at least one complex root. This includes polynomials with real coefficients as well. As a result, we can factor any polynomial $f(x)$ with complex coefficients in $\mathbb{C}$ into linear factors to obtain

$$f(x) = (x - \alpha_1)^{m_1}(x - \alpha_2)^{m_2} \cdots (x - \alpha_p)^{m_p}$$

where $\alpha_1, \alpha_2, \ldots, \alpha_p$ are in $\mathbb{C}$ and $m_i$ is the multiplicity of $\alpha_i$.

This result along with the Cayley-Hamilton Theorem is the key tool in producing the main result in this article. Cayley-Hamilton Theorem allows us to write a polynomial equation as a matrix equation. It states that [18]:

**Theorem 2.5** (Cayley-Hamilton Theorem). *A matrix satisfies its own characteristic equation.*

As a result, we may factor a matrix polynomial equation into linear factors and then leads to the product of matrices called $\kappa$-matrices. We have the following main definition of this work.

**Definition 2.6** ($\kappa$- matrix). Let $A$ be an $n \times n$ matrix and $\lambda \in \sigma(A)$. Then we call the matrix given by $A - \lambda I$ is the *eigenmatrix* or the *characteristic matrix* of $A$ and is denoted by $\kappa_\lambda(A)$. Due to this notation, we sometimes call the characteristic matrix a $\kappa$-matrix, as well.

The determinant of the $\kappa$-matrix is called the characteristic polynomial of $A$ and is denoted by $p_A(\lambda)$ or simply by $p(\lambda)$, if there is no ambiguity.

The matrix product is not commutative in general. However, when the matrices are of the form $A - \lambda I$, we have the following result.

**Lemma 2.7.** *The product of $\kappa$-matrices is commutative for a given matrix $A$.*

To establish the result it is sufficient to note that $(A - \lambda I)(A - \mu I) = A^2 - (\lambda + \mu)A + \lambda\mu I$ for a given matrix $A$ and two numbers $\lambda$ and $\mu$, real or complex.

The following result is fundamental in linear algebra and can be found, for example, in [31].

**Theorem 2.8.** *If $\mathbf{v}_1, \ldots, \mathbf{v}_r$ are the eigenvectors associated with the distinct eigenvalues $\lambda_1, \ldots, \lambda_r$ of an $n \times n$ matrix $A$, then they are linearly independent.*

*Proof.* The proof may be done using a contradiction argument. So, to the contrary, suppose the vectors $\mathbf{v}_1, \ldots, \mathbf{v}_r$ are linearly dependent. Since $\mathbf{v}_1$ is nonzero, and the set of vectors is finite, we may assume that one of the vectors is a linear combination of the preceding vectors. Let $p$ be the least index of such a vector, $\mathbf{v}_{p+1}$. Then there are $c_i$ such that

$$c_1\mathbf{v}_1 + \cdots + c_p\mathbf{v}_p = \mathbf{v}_{p+1} \tag{2.6}$$

Multiplying both sides of (2.6) by $A$ and using the fact that $A\mathbf{v}_i = \lambda_i\mathbf{v}_i$ for each $i$, we have

$$c_1\lambda_1\mathbf{v}_1 + \cdots + c_p\lambda_p\mathbf{v}_p = \lambda_{p+1}\mathbf{v}_{p+1} \tag{2.7}$$

Multiplying both sides of (2.6) by $\lambda_{p+1}$ and subtracting from (2.7), we have

$$c_1(\lambda_1 - \lambda_{p+1})\mathbf{v}_1 + \cdots + c_p(\lambda_p - \lambda_{p+1})\mathbf{v}_p = \mathbf{0} \tag{2.8}$$

Since the set $\{\mathbf{v}_1, \cdots, \mathbf{v}_p\}$ is linearly independent, this then implies that $\lambda_i - \lambda_{p+1}$ are all zero. But non of these factors are zero since the eigenvalues are taken to be distinct. Hence $c_i = 0$ for $i = 1, \ldots, p$. But then (2.6), says that $\mathbf{v}_{p+1} = \mathbf{0}$, which is impossible. Therefore the set $\{\mathbf{v}_1, \cdots, \mathbf{v}_r\}$ has to be linearly independent. This completes the proof. □

The following result is standard in the theory of linear algebra and can be found, for example, in [18, Section 5.2]. Not all the matrices are diagonalizable. Doing so amounts to find certain eigenvectors that form the transformation matrix $P$. The Theorem 2.8 says that the eigenvectors corresponding to different eigenvalues are linearly independent. But, is that enough to make a matrix diagonal. The next theorem establishes the requirements for a matrix to be diagonalizable.

**Theorem 2.9** (The Diagonalization Theorem)**.** *An $n \times n$ matrix $A$ is diagonalizable if and only if*

a) *$A$ has $n$ linearly independent eigenvectors.*
b) *The sum of the dimensions of the eigenspaces equals $n$. This happens if and only if (i) the characteristic polynomial factors completely into linear factors and (ii) the dimension of the eigenspace for each $\lambda_k$ equals the multiplicity of $\lambda_k$. In other words, the geometric multiplicity of each eigenvalue equals its algebraic multiplicity.*

*In general, the dimension of the eigenspace for $\lambda_k$, also called the geometric multiplicity, is less than or equal to the multiplicity (algebraic) of the eigenvalue $\lambda_k$.*

As a result of Theorems 2.8 and 2.9, if $A$ has $n$ distinct eigenvalues, then $A$ is diagonalizable. This is only a sufficient condition for $A$ to be diagonalizable.

With the help of this terminology, we have the following interesting result. That is, if a diagonalizable matrix has a spectrum of size 2, that is $|\sigma(A)| = 2$, we do not need to find the eigenvectors the usual way thanks to the following result, another version of which is called the *2-spectrum lemma* later.

**Theorem 2.10.** *Eigenvectors of a diagonalizable matrix $A$ with $|\sigma(A)| \leq 2$ are the nonzero columns of the $\kappa$-matrices, in reverse order (that is, of the other eigenvalue).*

Before producing a formal proof of Theorem 2.10, we will look at several cases to motivate the argument.

**Example 2.11.** As the first example, let $A = \begin{pmatrix} 2 & 1 \\ 1 & 2 \end{pmatrix}$. This has two distinct eigenvalues, $\lambda_1 = 1$ and $\lambda_2 = 3$. Thus, according to Theorem 2.9, it is diagonalizable. By the definition of $\kappa$-matrices, we have

$$\kappa_{\lambda=1}(A) = A - 1 \cdot I = \begin{pmatrix} 1 & 1 \\ 1 & 1 \end{pmatrix} \quad \text{and} \quad \kappa_{\lambda=3}(A) = A - 3 \cdot I = \begin{pmatrix} -1 & 1 \\ 1 & -1 \end{pmatrix}$$

It is easy to see that the columns of $\kappa$-matrices are linearly dependent, if not the same. The first column of $\kappa_1(A)$ is $v_1 = \begin{bmatrix} 1 \\ 1 \end{bmatrix}$, which can be shown to be an eigenvector of $A$ for the eigenvalue $\lambda_2 = 3$. The second column of $\kappa_3(A)$ is $\begin{bmatrix} 1 \\ -1 \end{bmatrix}$, which is an associated eigenvector for $\lambda_1 = 1$.

One may argue that the method works only for symmetric matrices. However, it turns out that the approach also works for non-symmetric matrices and will be considered later in the sequel. The interested readers are asked to consider the matrix $\mathbf{A} = \left(\begin{smallmatrix} 3 & 1 \\ 2 & 4 \end{smallmatrix}\right)$ to see it for themselves. As the second example we consider a situation where there are repeated eigenvalues.

**Example 2.12.** Let $A = \begin{pmatrix} 3 & -1 \\ 1 & 1 \end{pmatrix}$. This has $\lambda = 2$ with multiplicity 2. Then,

$$\kappa_2(A) = A - 2I = \begin{pmatrix} 1 & -1 \\ 1 & -1 \end{pmatrix}.$$

The first column of $\kappa_2(A)$ is $v_1 = [1 \ \ 1]^T$, which can be shown to be an eigenvector of $A$. The second column of $\kappa_2(A)$ is $[-1 \ \ -1]^T$, which is a multiple of $v_1$, hence is another eigenvector of $A$. However, the eigenspace of $\lambda = 2$ is one dimensional, thus is spanned by only one vector. As a result, the matrix cannot be diagonalizable due to insufficient number of linearly independent eigenvectors.

It turns out that this result is valid for any 2-by-2 matrix with $|\sigma(A)| \leq 2$ and was proved in [29]. The proof is calculation based and utilizes the properties of quadratic formula.

**Theorem 2.13** (Diagonalizing $2 \times 2$ matrices)**.** *Let $A = \begin{pmatrix} a & b \\ c & d \end{pmatrix}$ with distinct eigenvalues $\lambda_1$ and $\lambda_2$. Then an eigenvector of $\lambda_1$ is $v_1 = \begin{bmatrix} a - \lambda_2 \\ c \end{bmatrix}$ and that of $\lambda_2$ is $v_2 = \begin{bmatrix} b \\ d - \lambda_1 \end{bmatrix}$. Then, we have*

$$P = \begin{pmatrix} a - \lambda_2 & b \\ c & d - \lambda_1 \end{pmatrix} \quad \text{and} \quad D = \begin{pmatrix} \lambda_1 & 0 \\ 0 & \lambda_2 \end{pmatrix}. \tag{2.9}$$

*In case of $v_1 = \mathbf{0}$, we use $v_1 = \begin{bmatrix} b \\ d - \lambda_2 \end{bmatrix}$. Similarly, we redefine $v_2$ in case of $v_2 = \mathbf{0}$. If there are no such nonzero $v_1 \neq kv_2$, $k \in \mathbb{C}$, the matrix $A$ is not diagonalizable.*

On a different point of view, the diagonal matrix $D = \begin{pmatrix} d & 0 \\ 0 & d \end{pmatrix}$ can be trivially diagonalized by using any invertible matrix $P$ of size $2 \times 2$ and by writing $D = PDP^{-1}$. Thus, the case of $|\sigma(A)| = 1$ is not very interesting. For $n > 2$, the case of non-diagonalizable matrix is complicated and leads to generalized eigenvectors and will be discussed in detail in another section.

It is an interesting question why certain matrices are not diagonalizable. It is due to the relationship among the columns of the corresponding $\kappa$-matrices. For example, if the complementary (other) $\kappa$-matrix has rank 1, then the eigenspace should be one dimensional. The rank of a matrix is a properties of its relationships among columns.

To our surprise, the idea of Theorem 2.13 works for $3 \times 3$ matrices with $|\sigma(A)| = 2$ as well. To see that consider the following.

**Example 2.14.** Let $A = \begin{pmatrix} 0 & -1 & 1 \\ -2 & -1 & 2 \\ -1 & -1 & 2 \end{pmatrix}$. The characteristic equation of $A$ is $\lambda^3 - \lambda^2 - \lambda + 1 = 0$ that gives the eigenvalues $\lambda_1 = 1$ with multiplicity 2 and $\lambda_2 = -1$. The corresponding eigenmatrices are:

$$\overset{\lambda_1 = 1}{\kappa_{\lambda_1}(A) = \begin{bmatrix} -1 & -1 & 1 \\ -2 & -2 & 2 \\ -1 & -1 & 1 \\ v_1 & v_2 & v_3 \end{bmatrix}} \quad \text{and} \quad \overset{\lambda_2 = -1}{\kappa_{\lambda_2}(A) = \begin{bmatrix} 1 & -1 & 1 \\ -2 & 0 & 2 \\ -1 & -1 & 3 \\ u_1 & u_2 & u_3 \end{bmatrix}}$$

It can be seen that the vectors $v_1, v_2$, and $v_3$ are linearly dependent and $v_1 = v_2 = -v_3$. The column space of $\kappa_1(A) = \{v | v = tv_3, \ t \in \mathbb{R}\}$. Thus it is of one dimensional. On the other hand, $u_2 = -\frac{1}{2}(u_1 + u_3)$ and $\text{span}(u_1, u_2, u_3) = \text{span}(u_2, u_3)$. Also, the vectors $u_2$ and $u_3$ are linearly independent. Thus, the column space of $\kappa_{-1}(A) = \text{span}(u_2, u_2)$ and the rank of $\kappa_{-1}(A)$ is 2. In addition, in Eq.(2.10) below, $w = \frac{1}{2}(u_3 - u_2)$.

While the diagonalization of a matrix is not unique, a possible diagonalization of $A$ may be given by

$$A = \begin{pmatrix} 0 & -1 & 1 \\ -2 & -1 & 2 \\ -1 & -1 & 2 \end{pmatrix} = \overset{P}{\begin{pmatrix} 1 & 1 & 1 \\ 2 & 1 & 2 \\ 3 & 2 & 1 \end{pmatrix}}_{u_3 \; w \; v_3} \overset{D}{\begin{pmatrix} 1 & 0 & 0 \\ 0 & 1 & 0 \\ 0 & 0 & -1 \end{pmatrix}}_{\lambda_1 \; \lambda_1 \; \lambda_2} \begin{pmatrix} -\frac{3}{2} & \frac{1}{2} & \frac{1}{2} \\ 2 & -1 & 0 \\ \frac{1}{2} & \frac{1}{2} & -\frac{1}{2} \end{pmatrix}. \tag{2.10}$$

*Remark* 2.15. According to the discuss above, it can be observed that, for a diagonalizable $3 \times 3$ matrix with spectrum of size 2, rank of one $\kappa$-matrix is one, while the other is of rank 2. In this case, any nonzero column of rank one matrix produces the eigenvector of $\lambda_2$ and any two nonzero columns of rank two eigenmatrix produces the eigenvectors of $\lambda_1$, which has algebraic multiplicity 2. Therefore, in $P$ above, the eigenvector corresponds to $\lambda_2$ is any nonzero multiple of any linearly independent column of $\kappa_{\lambda_1}(A)$ and the eigenvectors corresponds to $\lambda_1$ are any different linear combinations of any two linearly independent columns of $\kappa_{\lambda_2}(A)$.

Thus, in practice, we may look for simple and easy to work with eigenvectors for further calculations if necessary. The ultimate goal of finding eigenvectors would be to use in certain applications. Having something easy to work with will definitely provide added benefits in computational projects, specially in academia. We may introduce these techniques as classroom tools or "shortcut" methods for pedagogical purposes.

How about a matrix with some not so nice numbers! To that end, consider a matrix with some entries that are not integers. Let's study the following $3 \times 3$ matrix which has some rational entries.

**Example 2.16.** Let $B = \begin{pmatrix} \frac{3}{2} & -\frac{1}{2} & \frac{1}{2} \\ -1 & 1 & 1 \\ -\frac{1}{2} & -\frac{1}{2} & \frac{5}{2} \end{pmatrix}$, which has eigenvalues $\lambda_1 = 2$ with multiplicity 2 and $\lambda_2 = 1$. A possible diagonalization of $B$ would be

$$B = \begin{pmatrix} \frac{3}{2} & -\frac{1}{2} & \frac{1}{2} \\ -1 & 1 & 1 \\ -\frac{1}{2} & -\frac{1}{2} & \frac{5}{2} \end{pmatrix} = \overset{P}{\begin{pmatrix} 1 & 1 & 1 \\ 2 & 0 & 2 \\ 1 & 1 & 3 \end{pmatrix}}_{v_1 \; u_2 \; u_3} \overset{D}{\begin{pmatrix} 1 & 0 & 0 \\ 0 & 2 & 0 \\ 0 & 0 & 2 \end{pmatrix}}_{\lambda_2 \; \lambda_1 \; \lambda_1} \begin{pmatrix} \frac{1}{2} & \frac{1}{2} & -\frac{1}{2} \\ 1 & -\frac{1}{2} & 0 \\ -\frac{1}{2} & 0 & \frac{1}{2} \end{pmatrix}. \tag{2.11}$$

The corresponding $\kappa$-metrices are,

$$\kappa_{\lambda_1}(B) = \overset{\lambda_1 = 2}{\begin{bmatrix} -\frac{1}{2} & -\frac{1}{2} & \frac{1}{2} \\ -1 & -1 & 1 \\ -\frac{1}{2} & -\frac{1}{2} & \frac{1}{2} \end{bmatrix}}_{v_1 \; v_2 \; v_3} \quad \text{and} \quad \kappa_{\lambda_2}(B) = \overset{\lambda_2 = 1}{\begin{bmatrix} \frac{1}{2} & -\frac{1}{2} & \frac{1}{2} \\ -1 & 0 & 1 \\ -\frac{1}{2} & -\frac{1}{2} & \frac{3}{2} \end{bmatrix}}_{u_1 \; u_2 \; u_3}$$

This example has taken the maximum advantage of this approach just by taking multiples of the vectors appeared in $\kappa$-matrices for two eigenvalues. That is, $P$ is formed by the vectors $2v_3, 2u_2$ and $2u_3$, taken as columns that are corresponds to the different eigenvalues in $D$ in Eq. 2.11.

As the next example, let's turn our attention to complex eigenvalues and also possibly complex eigenvectors.

As complex eigenvalues of a real characteristic equation occur in conjugate pairs, there cannot be a third degree characteristic equation with three eigenvalues, two of which are real repeated while the other is complex or two complex repeated while the other is real. That is, there is no $3 \times 3$ real valued matrix with both real and complex eigenvalues, some of which has multiplicity 2. As for this reason, we consider a complex valued matrix with possibly real and complex eigenvalues, one of which is repeated.

**Example 2.17.** Let $C$ be a complex matrix with a diagonalization given by,

$$C=\begin{pmatrix}1-\frac{1}{2}i & \frac{1}{2}+i & -\frac{1}{2}-i\\ i & 3-i & -1\\ \frac{1}{2}+i & 1-\frac{1}{2}i & 1-\frac{1}{2}i\end{pmatrix}=\overset{P}{\underset{w_1\quad w_2\quad w_3}{\begin{pmatrix}1+2i & 1 & 1-i\\ 2 & 0 & 2\\ 2-i & i & 3+i\end{pmatrix}}}\overset{D}{\underset{\lambda_2\quad \lambda_1\quad \lambda_1}{\begin{pmatrix}1 & 0 & 0\\ 0 & 2-i & 0\\ 0 & 0 & 2-i\end{pmatrix}}}\overset{P^{-1}}{\begin{pmatrix}* & * & *\\ * & * & *\\ * & * & *\end{pmatrix}}$$

with the characteristic equation $-\lambda^3+5\lambda^2-7\lambda+3-i\left(2\lambda^2-6\lambda+4\right)=0$, which has eigenvalues $\lambda_1=2-i$ with multiplicity 2 and $\lambda_2=1$. The corresponding $\kappa$-matrices take the forms,

$$\kappa_{\lambda_1}(C)=\overset{\lambda_1=2-i}{\underset{v_1\qquad v_2\qquad v_3}{\begin{bmatrix}-1+\frac{1}{2}i & \frac{1}{2}+i & -\frac{1}{2}-i\\ i & 1 & -1\\ \frac{1}{2}+i & 1-\frac{1}{2}i & -1+\frac{1}{2}i\end{bmatrix}}}\quad\text{and}\quad \kappa_{\lambda_2}(C)=\overset{\lambda_2=1}{\underset{u_1\qquad u_2\qquad u_3}{\begin{bmatrix}-\frac{1}{2}i & \frac{1}{2}+i & -\frac{1}{2}-i\\ i & 2-i & -1\\ \frac{1}{2}+i & 1-\frac{1}{2}i & -\frac{1}{2}i\end{bmatrix}}}\tag{2.12}$$

It can be seen in (2.12) that $v_2=-v_3=v_1i$, thus the $\mathrm{col}(\kappa_{\lambda_1}(C))=\mathrm{span}(v2)$, taken over $\mathbb{C}$, and any nonzero multiple, real or complex, can be taken as the eigenvector of $\lambda_1$. For example, we have $w_1=2v_2$ in this example. Further, $\mathrm{col}(\kappa_{\lambda_2}(C))$ is spanned by any two nonzero columns of $\kappa_{\lambda_2}(C)$. Two eigenvectors of $\lambda_2=2-i$ are $w_2=\frac{1}{2}(u_1+u_3i)(1+i)$ and $w_3=\left(\frac{3}{2}-\frac{3}{2}i\right)u_1+\left(-\frac{1}{2}+\frac{3}{2}i\right)u_3$. Or, we could have simply taken $w_2=2u_1=[-i\ \ 2i\ \ 1+2i]^T$ and $w_3=-2u_3=[1+2i\ \ 2\ \ i]^T$ or any other two different linear combinations of $u_1$ and $u_3$. The standard row-reduction process of finding eigenvectors is involved and can be very time consuming, especially due to the involvement of complex numbers. Following this method we can simply avoid those unnecessary calculations, except finding the inverse of $P$.

Next, we look at the case of three distinct eigenvalues.

**Case 4: Three distinct eigenvalues.** Until now we have only considered at most two eigenvalues and we noticed that in each case the corresponding eigenvectors are the nonzero columns of the $\kappa$-matrices, in reverse order. Thus, the main observation is that $\kappa_{\lambda_1}(A)$ contains the eigenvectors of $\lambda_2$, while $\kappa_{\lambda_2}(A)$ contains that of $\lambda_1$. Under this case we study the case of three distinct eigenvalues of a $3\times3$ matrix and answer the question of how the observation in question carries over to three eigenvalues.

As the first example of this type, we consider the following.

**Example 2.18.** Let $E=\begin{pmatrix}4 & 0 & -1\\ 4 & 2 & -2\\ 5 & -1 & 0\end{pmatrix}$. The characteristic equation of $E$ is $-\lambda^3+6\lambda^2-11\lambda+6=0$ that gives the three eigenvalues $\lambda_1=1,\lambda_2=2$ and $\lambda_3=3$, each with multiplicity 1. The corresponding $\kappa$-matrices are given by

$$\kappa_{\lambda_1}(E)=\overset{\lambda_1=1}{\underset{v_1\quad v_2\quad v_3}{\begin{bmatrix}3 & 0 & -1\\ 4 & 1 & -2\\ 5 & -1 & -1\end{bmatrix}}}\quad \kappa_{\lambda_2}(E)=\overset{\lambda_2=2}{\underset{u_1\quad u_2\quad u_3}{\begin{bmatrix}2 & 0 & -1\\ 4 & 0 & -2\\ 5 & -1 & -2\end{bmatrix}}}\quad \kappa_{\lambda_3}(E)=\overset{\lambda_3=3}{\underset{w_1\quad w_2\quad w_3}{\begin{bmatrix}1 & 0 & -1\\ 4 & -1 & -2\\ 5 & -1 & -3\end{bmatrix}}}\tag{2.13}$$

According to (2.13), we can see that in each $\kappa$-matrix, only two columns are linearly independent, while the third is a linear combination of the other two. In particular, in this example we have $v_1=-2v_2-3v_3, u_1=-u_2-2u_3$ and $w_1=-2w_2-w_3$. Thus, to form $\mathrm{col}\kappa_{\lambda_i}(E)$, we may simply *pick any two nonzero columns* of each matrix.

In practice, we select columns with *least number of nonzero entries* and *least magnitude* with possible negation and scaling to make the calculation easier. In doing so, we may need to multiply entries by the *least common denominator* (LCD) of entries of the column to obtain integer entries. For example in this problem, we may pick $v_2,-v_3,-u_2,-u_3,-w_2$ and $-w_3$, as the best options for subsequent calculations.

As can be seen on Figure 4a, the column space of $\kappa_{\lambda_1}(E)$ contains the eigenvectors of $\lambda_2$ and $\lambda_3$. Similarly, $\mathrm{col}(\kappa_{\lambda_2}(E))$ contains the eigenvectors of $\lambda_1$ and $\lambda_3$. Thus, the intersection, $\mathrm{col}(\kappa_{\lambda_1}(E))\bigcap\mathrm{col}(\kappa_{\lambda_2}(E))$

contains the eigenspace of $\lambda_3$, which is a span of a single vector in this example. In particular, an eigenvector of $\lambda_3 = 3$ is $\mathrm{v}_3 = [1\ 2\ 1]^T$, which lies in the spaces of $\mathrm{col}(\kappa_{\lambda_1}(E))$ and $\mathrm{col}(\kappa_{\lambda_2}(E))$ since $\mathrm{v}_3 = -v_3 = u_2 - u_3$. A similar argument shows that $\mathrm{v}_1 = [1\ 2\ 3]^T$ is an eigenvector of $\lambda_1 = 1$ and $\mathrm{v}_2 = [1\ 1\ 2]^T$ is an eigenvector of $\lambda_2 = 2$. Hence a diagonalization of $E$ would be

$$E = \begin{pmatrix} 4 & 0 & -1 \\ 4 & 2 & -2 \\ 5 & -1 & 0 \end{pmatrix} = \overset{P}{\underset{\mathrm{v}_1 \quad \mathrm{v}_2 \quad \mathrm{v}_3}{\begin{pmatrix} 1 & 1 & 1 \\ 2 & 1 & 2 \\ 3 & 2 & 1 \end{pmatrix}}} \overset{D}{\underset{\lambda_1 \quad \lambda_2 \quad \lambda_3}{\begin{pmatrix} 1 & 0 & 0 \\ 0 & 2 & 0 \\ 0 & 0 & 3 \end{pmatrix}}} \begin{pmatrix} -\frac{3}{2} & \frac{1}{2} & \frac{1}{2} \\ 2 & -1 & 0 \\ \frac{1}{2} & \frac{1}{2} & -\frac{1}{2} \end{pmatrix}. \tag{2.14}$$

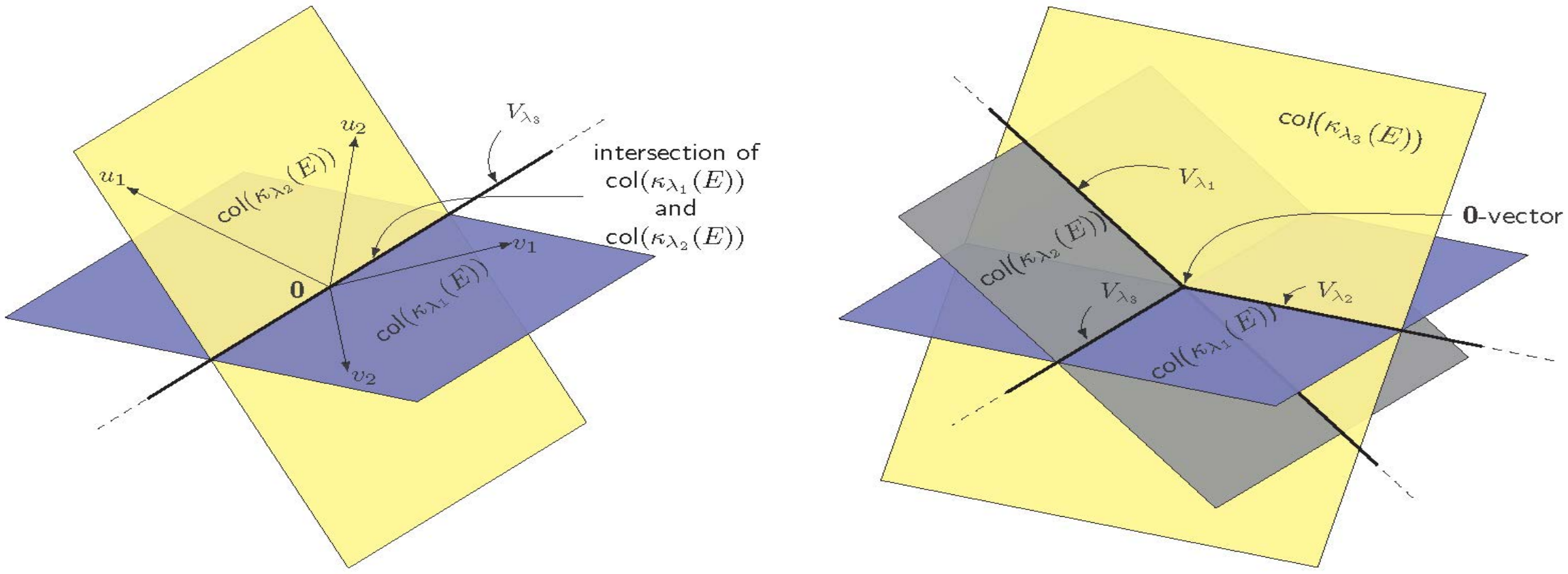

(A) Eigenspace of $\lambda_3$, $V_{\lambda_3}$ (B) Eigenspaces of $\lambda_1, \lambda_2$, and $\lambda_3$

FIGURE 4. Column spaces of $\kappa$-matrices

**Definition 2.19** (Complementary Eigenvalue)**.** Let $\{\lambda_1, \lambda_2, \ldots, \lambda_n\}$ be an indexed set, not necessarily distinct, of eigenvalues of an $n \times n$ matrix $A$. Then we call the eigenvalues $\lambda_i$ and $\lambda_j$ are *complementary* to one another if the index $i \neq j$.

According to this definition, an eigenvalue can be self-complementary when a matrix has repeated eigenvalues and will be useful when calculating eigenvectors for repeated eigenvalues.

With the help of above definition, the observation made in Example 2.18 can be summarized into the following process as a shortcut method for finding eigenvectors.

> **Quick Method of Finding Eigenvectors I**
>
> The eigenvectors can be found by a trail-and-error method by guessing the common vector(s) of the column spaces of $\kappa$-matrices of the complementary eigenvalues. For example, to find the eigenvector of $\lambda_1$, find a common vector of $\mathrm{col}(\kappa_{\lambda_2}(E))$ and $\mathrm{col}(\kappa_{\lambda_3}(E))$, if the matrix has three eigenvalues. The method can be easily extended to more than three eigenvalues.

### 2.1. Case 5: $4 \times 4$ matrix.

***Two eigenvalues.*** There are two cases to consider. If the matrix is diagonalizable, it should have either one eigenvalue with multiplicity one and another with multiplicity 3 or two eigenvalues each with multiplicity 2. Lets look at both of these cases in detail.

***Case 5(a)****: Multiplicity 1 and 3.* The eigenvalues may be taken as $\lambda_1$ and $\lambda_2$ with multiplicity 3. If the matrices are taken over $\mathbb{R}$, all the eigenvalues should be real since that occur as conjugate pairs for a real matrix.

**Example 2.20.** Let $F$ be given by

$$F = \begin{pmatrix} \frac{3}{2} & \frac{1}{2} & \frac{1}{2} & -\frac{1}{2} \\ 1 & 2 & 1 & -1 \\ \frac{3}{2} & \frac{3}{2} & \frac{5}{2} & -\frac{3}{2} \\ 2 & 2 & 2 & -1 \end{pmatrix}$$

Then $A$ has two eigenvalues $\lambda_1 = 2$ with multiplicity 1 and $\lambda_2 = 1$ with multiplicity 3. The corresponding $\kappa$-matrices are given by

$$\kappa_{\lambda_2}(F) = \overset{\lambda_2 = 1}{\begin{bmatrix} \frac{1}{2} & \frac{1}{2} & \frac{1}{2} & -\frac{1}{2} \\ 1 & 1 & 1 & -1 \\ \frac{3}{2} & \frac{3}{2} & \frac{3}{2} & -\frac{3}{2} \\ 2 & 2 & 2 & -2 \end{bmatrix}}_{\substack{v_1 \quad v_2 \quad v_3 \quad v_4}} \quad \text{and} \quad \kappa_{\lambda_1}(F) = \overset{\lambda_1 = 2}{\begin{bmatrix} -\frac{1}{2} & \frac{1}{2} & \frac{1}{2} & -\frac{1}{2} \\ 1 & 0 & 1 & -1 \\ \frac{3}{2} & \frac{3}{2} & \frac{1}{2} & -\frac{3}{2} \\ 2 & 2 & 2 & -3 \end{bmatrix}}_{\substack{u_1 \quad u_2 \quad u_3 \quad u_4}} \tag{2.15}$$

It can be seen that $v_1 = v_2 = v_3 = -v_4$, thus the $\text{col}(\kappa_{\lambda_2}(E)) = \text{span}(v_1)$. Therefore, an eigenvector of $\lambda_1 = 2$ may be taken as $\mathrm{v}_1 = v_1$ or $\mathrm{v}_1 = 2v_1$ . On the other hand, $u_4 = -\frac{1}{4}u_1 - \frac{1}{2}u_2 - \frac{3}{4}u_3$ and $u_1, u_2$ and $u_3$ are linearly independent. It can be shown that each of $u_i, i = 1, \ldots, 4$ is an eigenvector for $\lambda_1 = 2$. Thus, to form three linearly independent eigenvectors for $\lambda_2 = 1$, we may pick any three vectors from $\{u_1, u_2, u_3, u_4\}$ or any three nonzero linearly independent linear combinations of them. For example, we may select $v_2 = [1\ 0\ -1\ 0]^T = u_3 - u_1$, $v_3 = [1\ -1\ 0\ 0]^T = u_2 - u_1$ and $v_4 = [0\ 1\ 0\ 1]^T = -(u_2 + u_4)$ and form a diagonalization of $F$ as

$$F = \begin{pmatrix} \frac{3}{2} & \frac{1}{2} & \frac{1}{2} & -\frac{1}{2} \\ 1 & 2 & 1 & -1 \\ \frac{3}{2} & \frac{3}{2} & \frac{5}{2} & -\frac{3}{2} \\ 2 & 2 & 2 & -1 \end{pmatrix} = \overset{P}{\underset{\mathrm{v}_1 \quad \mathrm{v}_2 \quad \mathrm{v}_3 \quad \mathrm{v}_4}{\begin{pmatrix} 1 & 1 & 1 & 0 \\ 2 & 0 & -1 & 1 \\ 3 & -1 & 0 & 0 \\ 4 & 0 & 0 & 1 \end{pmatrix}}} \overset{D}{\underset{\lambda_1 \quad \lambda_2 \quad \lambda_2 \quad \lambda_2}{\begin{pmatrix} 2 & 0 & 0 & 0 \\ 0 & 1 & 0 & 0 \\ 0 & 0 & 1 & 0 \\ 0 & 0 & 0 & 1 \end{pmatrix}}} \begin{pmatrix} \frac{1}{2} & \frac{1}{2} & \frac{1}{2} & -\frac{1}{2} \\ \frac{3}{2} & \frac{3}{2} & \frac{1}{2} & -\frac{3}{2} \\ -1 & -2 & -1 & 2 \\ -2 & -2 & -2 & 3 \end{pmatrix}. \tag{2.16}$$

***Case 5(b)****: Multiplicity 2 and 2.* This case is much simpler than it appears since the eigenvectors are simply two nonzero nonparallel columns of $\kappa$-matrices. To see this consider the following.

**Example 2.21.** Let $G$ be given by

$$G = \begin{pmatrix} 1 & -\frac{1}{2} & \frac{1}{2} & -\frac{1}{2} \\ 0 & \frac{3}{2} & -\frac{1}{2} & \frac{1}{2} \\ -1 & -1 & 2 & 0 \\ -1 & -\frac{1}{2} & \frac{1}{2} & \frac{3}{2} \end{pmatrix}$$

Then $A$ has two eigenvalues $\lambda_1 = 1$ and $\lambda_2 = 2$ each with multiplicity 2. The corresponding $\kappa$-matrices are given by

$$\kappa_{\lambda_1}(G) = \overset{\lambda_1 = 1}{\begin{bmatrix} 0 & -\frac{1}{2} & \frac{1}{2} & -\frac{1}{2} \\ 0 & \frac{1}{2} & -\frac{1}{2} & \frac{1}{2} \\ -1 & -1 & 1 & 0 \\ -1 & -\frac{1}{2} & \frac{1}{2} & \frac{1}{2} \end{bmatrix}}_{\substack{v_1 \quad v_2 \quad v_3 \quad v_4}} \quad \text{and} \quad \kappa_{\lambda_2}(G) = \overset{\lambda_2 = 2}{\begin{bmatrix} -1 & -\frac{1}{2} & \frac{1}{2} & -\frac{1}{2} \\ 0 & -\frac{1}{2} & -\frac{1}{2} & \frac{1}{2} \\ -1 & -1 & 0 & 0 \\ -1 & -\frac{1}{2} & \frac{1}{2} & -\frac{1}{2} \end{bmatrix}}_{\substack{u_1 \quad u_2 \quad u_3 \quad u_4}} \tag{2.17}$$

It is easy to see that there are only *two* linearly independent columns in each $\kappa$-matrix and any two nonzero nonparallel columns of it can be taken as eigenvectors of the complimentary eigenvalues. For calculation purposes, suitable linear combinations of those two vectors can be selected. In particular, two eigenvalues of $\lambda_2 = 2$ can be taken as $\mathrm{v}_1 = -v_1 = [0\ 0\ 1\ 1]^T$ and $\mathrm{v}_2 = 2v_4 = [-1\ 1\ 0\ 1]^T$, while that of $\lambda_2 = 2$ is $\mathrm{v}_3 = -v_1 = [1\ 0\ 1\ 1]^T$ and $\mathrm{v}_4 = 2v_3 = [1\ -1\ 0\ 1]^T$. It is also interesting to see that two vectors are negative of each other, leaving only three different vectors in each characteristic matrix. Any two of those three nonzero vectors would work as eigenvectors.

***Three eigenvalues***. There is only one possibility with one eigenvalue repeated once along with two other eigenvalues. Two eigenvalues that are of multiplicity one may be complex conjugate in the case where the matrix is real.

***Case 5(c)***: *Multiplicity 2 and Multiplicity 1.*

**Example 2.22.** Let $H$ be given by

$$H = \begin{pmatrix} 3 & 2 & 5 & -5 \\ 3 & 4 & 7 & -7 \\ 4 & 4 & 10 & -9 \\ 6 & 6 & 14 & -13 \end{pmatrix}$$

$H$ has three eigenvalues $\lambda_1 = 1$ with multiplicity 2, $\lambda_2 = 0$ and $\lambda_3 = 2$ each with multiplicity 1. The corresponding eigenmatrices are given by

$$\kappa_{\lambda_1}(H) = \overset{\lambda_1 = 1}{\underset{\begin{matrix} u_1 & u_2 & u_3 & u_4 \end{matrix}}{\begin{bmatrix} 2 & 2 & 5 & -5 \\ 3 & 3 & 7 & -7 \\ 4 & 4 & 9 & -9 \\ 6 & 6 & 14 & -14 \end{bmatrix}}} \quad \kappa_{\lambda_2}(H) = \overset{\lambda_2 = 0}{\underset{\begin{matrix} v_1 & v_2 & v_3 & v_4 \end{matrix}}{\begin{bmatrix} 3 & 2 & 5 & -5 \\ 3 & 4 & 7 & -7 \\ 4 & 4 & 10 & -9 \\ 6 & 6 & 14 & -13 \end{bmatrix}}} \quad \kappa_{\lambda_3}(H) = \overset{\lambda_3 = 2}{\underset{\begin{matrix} w_1 & w_2 & w_3 & w_4 \end{matrix}}{\begin{bmatrix} 1 & 2 & 5 & -5 \\ 3 & 2 & 7 & -7 \\ 4 & 4 & 8 & -9 \\ 6 & 6 & 14 & -15 \end{bmatrix}}} \tag{2.18}$$

It can be seen that $\kappa_{\lambda_1}(H)$ has only *two* linearly independent columns, while $\kappa_{\lambda_2}(H)$ and $\kappa_{\lambda_3}(H)$ each has only *three* linearly independent columns. As mentioned several times earlier, for further calculation purposes, any two nonparallel vectors can be selected from the former and any three linearly independent vectors or linear combinations from the latter possibly with scaling and least number of nonzero terms. In Equation 2.18, in particular, we may select $u_1, u_3, v_1, v_2, v_3, w_1, w_2$ and $w_3$ or appropriate linear combinations. A new set of linear combinations such as $\mathrm{u}_1 = u_1, \mathrm{u}_2 = u_3 - 2u_1 = [1\ 1\ 1\ 2]^T, \mathrm{v}_1 = v_1, \mathrm{v}_2 = v_2 - v_1 = [-1\ 1\ 0\ 0]^T$ and $\mathrm{v}_3 = v_3 - 2v_1 = [-1\ 1\ 2\ 2]^T$ will be used for further calculations. The rest of the calculation goes as follows.

Since $\mathrm{col}(\kappa_{\lambda_1}(H))$ is spanned by the two vectors $\mathrm{u}_1$ and $\mathrm{u}_2$ and it contains the eigenvectors of the other two eigenvalues $\lambda_2 = 0$ and $\lambda_3 = 2$, they can be expressed as $v = a\mathrm{u}_1 + b\mathrm{u}_2$ for some $a, b \in \mathbb{R}$. The goal is to find such $a$ and $b$ so that it produces the eigenvectors for $\lambda_2$ and $\lambda_3$. To that end, first note that

$$v = a\mathrm{u}_1 + b\mathrm{u}_2 = a\begin{pmatrix} 2 \\ 3 \\ 4 \\ 6 \end{pmatrix} + b\begin{pmatrix} 1 \\ 1 \\ 1 \\ 2 \end{pmatrix} = \begin{pmatrix} 2a+b \\ 3a+b \\ 4a+b \\ 6a+2b \end{pmatrix} \tag{2.19}$$

Now, we use the fact that $v$ is an eigenvector of $H$ and $Hv = \lambda_i v$ for the two eigenvalues. Since we need only one relation between $a$ and $b$, we use only one row of $H$ to produce a condition on the two variables. We reach the same relation irrespective of the row that has been used.

- **For** $\lambda_2 = 0$: Using the first row of $H$ we have

  $$3(2a+b) + 2(3a+b) + 5(4a+b) - 5(6a+2b) = 0 \Longrightarrow a = 0.$$

  This also implies that $b$ is any real number. Thus, an eigenvector of $\lambda_2 = 0$ is $v = [1\ 1\ 1\ 2]^T$.
- **For** $\lambda_3 = 2$: Using the first row of $H$ again we have

  $$3(2a+b) + 2(3a+b) + 5(4a+b) - 5(6a+2b) = 2(3a+b) \Longrightarrow a = -b, b \in \mathbb{R}.$$

  In particular, taking $a = 1$ and $b = -1$, an eigenvector of $\lambda_3 = 2$ is $v = [1\ 2\ 3\ 4]^T$.

To find the eigenvectors of $\lambda_1$, we use only one of $\kappa$-matrices. The solution does not depend of the matrix used and we tend to use the matrix with entries close to zero to make the calculation easier. Using a similar argument and the redefined vectors, we may write the form of the eigenvectors of $\lambda_1$ (and also of $\lambda_3$, which has already been found) as

$$v = a\mathrm{v}_1 + b\mathrm{v}_2 + c\mathrm{v}_3 = a\begin{pmatrix}3\\3\\4\\6\end{pmatrix} + b\begin{pmatrix}-1\\1\\0\\0\end{pmatrix} + c\begin{pmatrix}-1\\1\\2\\2\end{pmatrix} = \begin{pmatrix}3a-b-c\\3a+b+c\\4a+2c\\6a+2c\end{pmatrix}. \tag{2.20}$$

Now, for $\lambda_1 = 1$, we have, due to $Hv = \lambda_1 v$ and the first row of $H$ that

$$3(3a-b-c)+2(3a+b+c)+5(4a+2c)-5(6a+2c) = 3a-b-c \Longrightarrow a = 0.$$

This also implies that $b$ and $c$ can be any real numbers. Taking $b = 0$ and $c = 0$ separately, we may produce two eigenvectors of $\lambda_1$. Thus, we have $v = [-1\ 1\ 0\ 0]^T$ and $v = [-1\ 1\ 2\ 2]^T$ as two eigenvalues. These choices of eigenvectors lead to the diagonalization given by

$$H = \begin{pmatrix}1&1&-1&-1\\1&2&1&1\\1&3&0&2\\2&4&0&2\end{pmatrix}\begin{pmatrix}0&0&0&0\\0&2&0&0\\0&0&1&0\\0&0&0&1\end{pmatrix}\begin{pmatrix}-1&-1&-3&3\\1&1&2&-2\\0&1&0&-\frac{1}{2}\\-1&-1&-1&\frac{3}{2}\end{pmatrix}.$$

Other options may include $\mathrm{v}_3 = \frac{1}{2}(v_3 - 2v_1 - \mathrm{v}_2) = [0\ 0\ 1\ 1]^T$.

*Remark* 2.23. It is worth noting that all the eigenvectors are produced by only two $\kappa$-matrices. In theory, one $\kappa$-matrix contains all the associated eigenvectors as linear combination of columns, except for the eigenvalue that generates the $\kappa$-matrix itself, a fact we prove in the sequel.

***Case 5(d)****: Four eigenvalues: Each Multiplicity 1.* As the final observation, even in the case of eigenvalues with multiplicity one, we may do the same as we have done earlier by taking complementary eigenvectors as linear combinations of columns of the $\kappa$-matrices. We consider a detailed problem in the illustrative example section later.

*Remark* 2.24. Since any nonzero multiple of an eigenvector is still an eigenvector, we may scale them before using in $P$ to make the calculation simpler. We will discuss that in an example later. As another note, we do not distinguish eigenvalues as real or complex in Theorem **??** and the result is also valid for complex eigenvalues as well. As the third observation, we notice that the eigenvalues in $D$ appears in *reverse order.*

## 3. Main Results

The following result is elementary and can be proved directly using the characteristic equation of the matrix.

**Theorem 3.1** (Shifting Theorem)**.** *If a square-matrix $A$ has an eigenvalue $\lambda$, then the matrix $A - \mu I$ has an eigenvalue $\lambda - \mu$.*

*Proof.* This can be seen since the characteristic polynomial of $A$ is $P_\lambda = \lambda^n + a_{n-1}\lambda^{n-1} + \cdots + a_0$ and that of $A - \mu I$ is $P_{\lambda-\mu} = (\lambda-\mu)^n + a_{n-1}(\lambda-\mu)^{n-1} + \cdots + a_0$ with the same coefficients in both polynomials. □

One of the major properties of an eigenvalue is that when subtracted it from the diagonal entries of a matrix, it produces a non-invertible matrix, a key property to find eigenvectors. In that line of argument, we have the following result.

**Theorem 3.2.** *$\kappa_\lambda$-matrices are singular if $\lambda$ is an eigenvalue.*

This is clear since any eigenvalue is a solution of the corresponding characteristic equation. $\kappa_\lambda$-matrices always have 0 as an eigenvalue. The other eigenvalues can be found from those of $A$ using *the Shifting Theorem* above.

As a consequence of Theorem 3.2, $\kappa_\lambda$-matrices of a given matrix are always singular. Thus, the columns of each $\kappa_\lambda$-matrix are linearly dependent. Therefore, there is at least one column that can be written as a linear combination of rest of the columns. This then leads to fact that 0 can be written as a nontrivial linear combination of columns of $\kappa_\lambda$-matrix. The weights produce an eigenvector of $\lambda$. So, we have the result.

**Theorem 3.3.** *$\kappa$-matrices has 0 as an eigenvalue.*

### 3.1. Eigenvectors of $A$ and $A - \lambda I$.

We have been overlooking the use of columns to represent eigenvectors for centuries. The standard approach is to use the Gaussian elimination to reduce the system $(A - \lambda I)v = 0$ to its echelon-form and then solve for $v$. The Gaussian elimination is an operation on rows rather than on columns as we use in finding the eigenvectors.

3.1.1. *What happens when* 0 *is an eigenvalue.* When 0 is an eigenvalue, the corresponding matrix $A$ is singular, meaning that the determinant is 0. In other words, at least one column of $A$ is linearly depend on other columns.

**Theorem 3.4.** *If an $n \times n$ matrix $A$ has only two eigenvalues, $\lambda_1$ of multiplicity $m$ and $\lambda_2$ of multiplicity $n - m$, then the matrix $\kappa_{\lambda_2}(A)$ has* 0 *as an eigenvalue with algebraic multiplicity $m$ and $\kappa_{\lambda_1}(A)$ has* 0 *as an eigenvalue with algebraic multiplicity $n - m$.*

This is the reason why these $\kappa$-matrices have only $m$ and $n - m$ many linearly independent columns, respectively producing the eigenvectors for the associated complementary eigenvalues.

**Lemma 3.5** (2 – Spectrum Lemma I)**.** *Suppose an $n \times n$ matrix $A$ has a spectrum of 2, with $\lambda_1$ and $\lambda_2$ being the two eigenvalues. Then, $A$ is diagonalizable if and only if*

$$(A - \lambda_1 I) \cdot (A - \lambda_2 I) = 0. \tag{3.1}$$

*In addition, the nonzero columns of $A - \lambda_2 I$ are the eigenvectors of $A$ for the eigenvalue $\lambda_1$ and vice verse.*

The result can also be expressed in a more compact way as $\kappa_{\lambda_1}(A) \cdot \kappa_{\lambda_2}(A) = 0_n$ or $\langle \kappa_{\lambda_1}(A), \kappa_{\lambda_2}(A) \rangle = \mathrm{tr}((\kappa_{\lambda_1}(A))^T \kappa_{\lambda_1}(A))$, in the classical sense.

*Proof.* Let $A$ be diagonalizable. Then due to the Diagonalization Theorem, there is an eigenbasis for the underlying vector space. Let it be $\mathcal{B} = \{\mathrm{v}_1, \ldots, \mathrm{v}_n\}$, where each of $\mathrm{v}_i$ is an eigenvector of one of the eigenvalues, $\lambda_1$ or $\lambda_2$. Due to the commutative property of the $\kappa$-matrices, $(A - \lambda_1 I)(A - \lambda_2 I) = (A - \lambda_2 I)(A - \lambda_1 I)$. Thus, $(A - \lambda_1 I)(A - \lambda_2 I)\mathrm{v}_i = 0$ for each $i = 1, \ldots, n$, by considering the appropriate order of operations as needed. Now, let $v$ be any vector in the underlying vector space. Since $\mathcal{B}$ is a basis, $v$ can be written as a linear combination as in $v = c_1\mathrm{v}_1 + \cdots + c_n\mathrm{v}_n$. Then we have

$$\begin{aligned}(A - \lambda_1 I)(A - \lambda_2 I)v &= M\left(c_1\mathrm{v}_1 + \cdots + c_n\mathrm{v}_n\right), \text{ where } M = (A - \lambda_1 I)(A - \lambda_2 I) \\ &= c_1 M\mathrm{v}_1 + \cdots + c_n M\mathrm{v}_n, \\ &= 0.\end{aligned}$$

Since this is true for any $v$, we arrive at the desired result. To show that the nonzero columns of $(A - \lambda_2 I)$ are eigenvectors associated with the eigenvalue $\lambda_1$, we let $(A - \lambda_2 I) = (v_1, \ldots, v_n)$. Then, we have $(A - \lambda_1 I)(v_1, \ldots, v_n) = 0$, which then implies that $(A - \lambda_1 I)v_i = 0$ for each $i = 1, \ldots, n$. Therefore, $v_i$ is either a zero vector or an eigenvector of $A$ for the eigenvalue $\lambda_1$. Reversing the order in $(A - \lambda_1 I)(A - \lambda_2 I) = 0$, we can similarly prove that the nonzero columns of $(A - \lambda_1 I)$ are the eigenvector of $\lambda_2$. The other direction of the proof is straightforward from the Rank Theorem [31, p. 235]. This completes the proof. □

Alternatively, we can express the same result above in another way as given below.

**Lemma 3.6** (2 – Spectrum Lemma II)**.** *If an $n \times n$ diagonalizable matrix $A$ has a spectrum of 2, then the nonzero columns of the $\kappa$-matrices are the eigenvectors of the complementary eigenvalue.*

*That is, the nonzero columns of $A - \lambda_1 I$ are the eigenvectors of $\lambda_2$ and that of $A - \lambda_2 I$ are the eigenvectors of $\lambda_1$, where $\lambda_1$ and $\lambda_2$ are distinct eigenvalues of $A$.*

The Lemma 3.5 has several interesting applications and are discussed in the next subsection. For pedagogical purposes and to see how the algebra works here, we consider several examples with direct calculation based approach, even though the general case is relatively easier to prove.

3.1.2. *Special Cases.* We consider several consequences of Lemma 3.5 for $2\times 2$ and higher order matrices. We may combine several $\kappa$-matrices into one matrix that would produce all the eigenvectors, if all the columns are nonzero. Another reason for this to work is that the number of linearly independent eigenvectors needed to form the diagonal matrix is the same as the dimension of the matrix. Thus, subtracting the complementary eigenvalues from the diagonal entries will produce the associated eigenvectors as columns, if they are nonzero. The method applies to any diagonalizable matrix of spectrum 2. Let's look at several cases in detail to further clarify the underlying argument.

***Case:*** $2\times 2$ *Matrices.* The case of $2\times 2$ has been proved in [29]. There is only one possibility with two eigenvalues. So, it is easy to produce eigenvectors by subtracting complementary eigenvalues from the diagonal entries, if the resulting columns are nonzero. For example, the eigenvectors of $\mathbf{A} = \begin{pmatrix} a & b \\ c & d \end{pmatrix}$ are the nonzero columns of

$$\kappa_{\lambda_1,\lambda_2}(A) = \begin{pmatrix} a-\lambda_2 & b \\ c & d-\lambda_1 \end{pmatrix}$$

if the corresponding eigenvalues are $\lambda_1$ and $\lambda_2$. In case of $A$ is diagonal with repeated eigenvalue $d\neq 0$, we subtract $\dfrac{d}{2}$ from the diagonal entries. This then give $[1\ 0]^T$ and $[0\ 1]^T$ as eigenvectors after scaling.

***Case:*** $3\times 3$ *Matrices.* A $3\times 3$ matrix has only one case of spectrum 2 with one eigenvalue being repeated. Suppose $A$ has two distinct eigenvalues $\lambda_1$ and $\lambda_2$ with latter being of multiplicity 2. Then the corresponding eigenvectors of $\lambda_1, \lambda_2$ and $\lambda_2$ are the nonzero columns of the matrix

$$\kappa_{\lambda_1,\lambda_2,\lambda_2}(A) = \begin{pmatrix} a_{11}-\lambda_2 & a_{12} & a_{13} \\ a_{21} & a_{22}-\lambda_1 & a_{23} \\ a_{31} & a_{32} & a_{33}-\lambda_1 \end{pmatrix},$$

respectively.

For pedagogical purposes we also prove the case for $3\times 3$ matrices by a purely computational argument. When the matrix is $3\times 3$, there is only one case to consider for $|\sigma(A)| = 2$. Without loss of generality, suppose the eigenvalues are $\lambda_1$ and $\lambda_2$ and $\lambda_2$ has multiplicity 2. Let the matrix $A$ be given by

$$A = \begin{pmatrix} a_{11} & a_{12} & a_{13} \\ a_{21} & a_{22} & a_{23} \\ a_{31} & a_{32} & a_{33} \end{pmatrix} \tag{3.2}$$

Then, using the properties of the characteristic polynomial of $A$, we have

$$\lambda_1 + 2\lambda_2 = \text{trace}(A) = a_{11} + a_{22} + a_{33} \tag{3.3}$$

$$2\lambda_1\lambda_2 + \lambda_2^2 = \sum_{1\leq i\leq j\leq 3} a_{ii}a_{jj} - a_{ij}a_{ji} \tag{3.4}$$

$$\lambda_1\lambda_2^2 = \det(A). \tag{3.5}$$

We need to show that the vector $\mathbf{v}_1 = [a_{11}-\lambda_2, a_{21}, a_{31}]^T$ is an eigenvector of $A$ for the eigenvalue $\lambda_1$. In other words, we need to show that the product $(A-\lambda_1 I)\cdot \mathbf{v}_1 = \mathbf{0}$. To that end, we prove that each component of the resulting vector is zero. There are two cases to consider. First, consider the first element of the vector. Then we have using (3.3),

$$\begin{aligned} \text{row 1 of } (A-\lambda_1 I)\begin{pmatrix} a_{11}-\lambda_2 \\ a_{21} \\ a_{31} \end{pmatrix} &= (a_{11}-\lambda_1)(a_{11}-\lambda_2) + a_{12}a_{21} + a_{13}a_{31} \\ &= (2\lambda_2 - a_{22} - a_{33})(a_{11}-\lambda_2) + a_{12}a_{21} + a_{13}a_{31} \\ &= [a_{13}a_{31} - (a_{11}-\lambda_2)(a_{33}-\lambda_2)] + [a_{12}a_{21} - (a_{11}-\lambda_2)(a_{22}-\lambda_2)] \\ &= -C_{22} - C_{33} = 0 + 0 = 0, \end{aligned} \tag{3.6}$$

where $C_{ii}$ is the cofactor of the $a_{ii}$ entry of $A$. Each cofactor is zero since the rank of $\kappa_{\lambda_2}(A) = 1$ because the eigenvalue $\lambda_2$ has a multiplicity 2. Similarly, we can show that the row 2 of $(A - \lambda_1 I)\,\mathbf{v}_1 = a_{23}a_{31} - a_{21}(a_{33} - \lambda_2) = -C_{12} = 0$ and row 3 of $(A - \lambda_1 I)\,\mathbf{v}_1 = a_{32}a_{21} - a_{31}(a_{22} - \lambda_2) = -C_{13} = 0$. Thus, row 2 of $(A - \lambda_1 I)\,\mathbf{v}_1 = \mathbf{0}$. Therefore $\mathbf{v}_1$ is an eigenvector of $A$ for the eigenvalue $\lambda_1$. Similarly, we can show that the vectors $[a_{12}, a_{22} - \lambda_2, a_{32}]^T$ and $[a_{13}, a_{23}, a_{33} - \lambda_2]^T$ are either zero or an eigenvector of $\lambda_1$.

Now, let $\mathbf{v}_2 = [a_{11} - \lambda_1, a_{21}, a_{31}]^T$ be a vector. We show that this is an eigenvector of $\lambda_2$. As in the previous case, we see that row 1 of $(A - \lambda_2 I)\,\mathbf{v}_2 = (a_{11} - \lambda_1)(a_{11} - \lambda_2) + a_{12}a_{21} + a_{13}a_{31} = 0$. The vanishing of other components of $\mathbf{v}_2$ can be proved similarly using the fact that the rank of $\kappa_{\lambda_2}(A) = 1$ in all these cases. Hence $\mathbf{v}_2$ is an eigenvector of $A$ for the eigenvalue $\lambda_2$. Other columns of $\kappa_{\lambda_2}(A)$ can be treated similarly.

The case of spectrum 3 will be much involved and leave as an interesting exercise for interested readers.

***Case:*** *4 × 4 Matrices.* There are two cases of 4 × 4 matrices with spectrum 2 to consider, one with multiplicities one and three and another with multiplicities two and two. Let the distinct eigenvalues of $A$ be $\lambda_1$ and $\lambda_2$ as before. The corresponding $\kappa$matrices for the two cases are given by

$$\kappa_{\lambda_1,\lambda_2,\lambda_2,\lambda_2,}(A) = \begin{pmatrix} a_{11}-\lambda_2 & a_{12} & a_{13} & a_{14} \\ a_{21} & a_{22}-\lambda_1 & a_{23} & a_{24} \\ a_{31} & a_{32} & a_{33}-\lambda_1 & a_{34} \\ a_{41} & a_{42} & a_{43} & a_{44}-\lambda_1 \end{pmatrix}$$

and

$$\kappa_{\lambda_1,\lambda_1,\lambda_2,\lambda_2,}(A) = \begin{pmatrix} a_{11}-\lambda_2 & a_{12} & a_{13} & a_{14} \\ a_{21} & a_{22}-\lambda_2 & a_{23} & a_{24} \\ a_{31} & a_{32} & a_{33}-\lambda_1 & a_{34} \\ a_{41} & a_{42} & a_{43} & a_{44}-\lambda_1 \end{pmatrix},$$

respectively. The (nonzero) columns are the eigenvectors of the associated complimentary eigenvalues. For example, the first column of $\kappa_{\lambda_1,\lambda_1,\lambda_2,\lambda_2,}(A)$ is an eigenvector of $\lambda_1$, which we obtained by subtracting the complementary eigenvalue $\lambda_2$ by the diagonal entry of that column.

The case of spectrum 3 and 4 for 4 × 4 matrices will be discussed in the next section. The method can be extended to any matrix of spectrum 2 irrespective of the dimension of the matrix.

3.2. **Relationship between col($A$) and col($A-\lambda I$).** When $\lambda$ is an eigenvalue of $A$, there is an interesting resemblance of eigenvectors of $A$ and the column spaces of $A - \lambda I$as given in the following main result of the paper. We call the theorem *The Eigenmatrix Theorem* or *The Characteristic Matrix Theorem* or *The $\kappa$-matrix Theorem* interchangeably.

**Lemma 3.7** (Eigenmatrix Lemma)**.** *If $\lambda$ is an eigenvalue of $A$, the eigenvector(s) of $A$ lie in the col($A - \lambda I$) for complementary eigenvalues.*

*That is, if $\mu \neq \lambda$ is another eigenvalue of $A$, the associated eigenvector of $\mu$ lies in col($A - \lambda I$).*

*Remark* 3.8. We will reproduce this result as the Corollary 3.11 later and provide a proof there. Interested readers are asked to provide a proof without using the Eigenmatrix Theorem below.

The following is the main result of the paper and can be called either *The Eigenmatrix Theorem* or the *$\kappa$-Matrix Theorem* for a reason apparent from the context.

**Theorem 3.9** (The Eigenmatrix Theorem ). *Let $A$ be an $n \times n$ matrix. Let $\lambda_i, 1 \le i \le p \le n$ be the distinct eigenvalues of $A$. Then, $A$ is diagonalizable if and only if*

$$(A - \lambda_1 I)\,(A - \lambda_2 I) \cdots (A - \lambda_p I) = 0 \tag{3.7}$$

*or*

$$\kappa_{\lambda_1}(A) \cdot \kappa_{\lambda_2}(A) \cdots \kappa_{\lambda_p}(A) = 0 \tag{3.8}$$

*In addition, the nonzero column vectors of the matrix*

$$\kappa_{\lambda_1}(A) \cdots \kappa_{\lambda_{i-1}}(A) \cdot \kappa_{\lambda_{i+1}}(A) \cdots \kappa_{\lambda_p}(A) \tag{3.9}$$

*are the (right) eigenvectors of $A$ for the eigenvalue $\lambda_i$ for $i = 2, \ldots, p-1$. The case for $i = 1$ and $p$ can be considered similarly.*

*Proof.* Let $p(\lambda) = \lambda^n + a_{n-1}\lambda^{n-1} + \cdots + a_0$ be the characteristic polynomial of $A$. Then due to the Cayley-Hamilton Theorem for matrices, we have

$$p(A) = A^n + a_{n-1}A^{n-1} + \cdots + a_0 I = 0. \tag{3.10}$$

Suppose $A$ has $n$ distinct eigenvalues $\lambda_i, i = 1, 2, \ldots, n$. Then by the Fundamental Theorem of Algebra, we may factor the left hand side of (3.10) as

$$(A - \lambda_1 I)\,(A - \lambda_2 I) \cdots (A - \lambda_p I) = 0. \tag{3.11}$$

This proves the claim for the case of $n$ distinct eigenvalues. Now suppose, some eigenvalues have algebraic multiplicities greater than 1. For simplicity, assume only $\lambda_p$ has a multiplicity greater than 1 since the other cases are similar due to the fact that the matrices are considered over commutative rings $\mathbb{R}$ or $\mathbb{C}$. Then the multiplicity of $\lambda_p$ should be $n - (p-1)$ and the Eq. 3.11 would then be

$$(A - \lambda_1 I)\,(A - \lambda_2 I) \cdots (A - \lambda_p I)^{n-p+1} = 0. \tag{3.12}$$

This can be factored into

$$[(A - \lambda_1 I)\,(A - \lambda_2 I) \cdots (A - \lambda_p I)]\,(A - \lambda_p I)^{n-p} = 0. \tag{3.13}$$

Unfortunately, this does not prove the claim that $(A - \lambda_1 I) \cdots (A - \lambda_p I) = 0$ since Eq. 3.13 is a matrix equation. However, it can be proved by following the discussion after Theorem 3 of Section 6.3 in [27]. We reproduce the proof here for completeness.

Now, suppose $A$ is diagonalizable and $\lambda_i, i = 1, 2, \ldots, p$ are the distinct eigenvalues of $A$. Then, due to Theorem 2.9, there is an eigenbasis for the underlying vector space. Let it be $\mathcal{B} = \{\mathbf{v}_i,\ \mathbf{v}_2,\ \ldots, \mathbf{v}_n\}$. Then it can be proved that $(A - \lambda_1 I)\,(A - \lambda_2 I) \cdots (A - \lambda_p I)\,\mathbf{v}_i = 0$ for any $i \in \{1, 2, \ldots, n\}$. To see this, consider

$$(A - \lambda_1 I)\,(A - \lambda_2 I) \cdots (A - \lambda_p I)\,\mathbf{v}_i \tag{3.14}$$

Without loss of generality, let $\lambda_i$ be the eigenvalue corresponding to $\mathbf{v}_i$. Since the multiplication of $(A - \lambda_i I)$ matrices in (3.14) is commutative, rearranging terms if necessary, we have

$$(A - \lambda_1 I)\,(A - \lambda_2 I) \cdots (A - \lambda_i I)\,\mathbf{v}_i = 0 \tag{3.15}$$

since $(A - \lambda_i I)\,\mathbf{v}_i = 0$. Similarly, we can show that similar annihilation hold for each $\mathbf{v}_i, i = 1, \ldots, n$. Now, to prove (3.7), we need to show that $(A - \lambda_1 I)\,(A - \lambda_2 I) \cdots (A - \lambda_p I)\,v = 0$ for any $v$ in the underlying vector space. To than end, since $\mathcal{B}$ is a basis, $v$ can be written as $v = c_1\mathbf{v}_1 + \cdots + c_n\mathbf{v}_n$. Then,

$$\begin{aligned}(A - \lambda_1 I)\,(A - \lambda_2 I) \cdots (A - \lambda_p I)\,v &= M(c_1\mathbf{v}_1 + \cdots + c_n\mathbf{v}_n)\\ &= c_1 M\mathbf{v}_1 + \cdots + c_n M\mathbf{v}_n\\ &= 0,\end{aligned} \tag{3.16}$$

since $M\mathbf{v}_i = 0, i = 1, \ldots, n$, where $M = (A - \lambda_1 I)\,(A - \lambda_2 I) \cdots (A - \lambda_p I)$. Eq. 3.16, then implies that

$$(A - \lambda_1 I)\,(A - \lambda_2 I) \cdots (A - \lambda_p I) = 0$$

This establishes (3.7).

Now, to prove (3.9), without loss of generality, we prove that $(A-\lambda_1 I)u_i = 0$, where $u_i$ is any nonzero column of the matrix $(A-\lambda_2 I)\cdots(A-\lambda_p I)$. We can do this since the matrices are considered over a commutative rings. Then we have by Eq. 3.7

$$(A-\lambda_1 I)\, u_i = 0.$$

Since $u_i \neq 0$, this shows that $u_i$ is an associated eigenvector of $A$ for $\lambda_1$. The sufficiency is straightforward and left for the reader to prove. □

It is worth noting that the diagonalizability of $A$ is needed in the first part of the proof to write any vector that is in the underlying space as a linear combination of eigenbasis vectors. For this reason we give the following definition.

**Definition 3.10.** A diagonalizable matrix is called *complete* and an nondiagonalizable matrix is referred to as *defective.*

The $v_i$ discussed above has a further explanation. It can be shown that such a eigenvector is in the column space of all complementary eigenvalues. We have the following result.

**Corollary 3.11.** *Let $v_i$ be an eigenvector of $A$ for the eigenvalue $\lambda_i$. Then, $v_i \in col\big(\kappa_{\lambda_j}(A)\big)$ for $j \neq i$.*

*Proof.* Let $\lambda_j \neq \lambda_i$ and $\kappa_{\lambda_j}(A)$ be the eigenmatrix of $\lambda_j$. Since the multiplication of $\kappa$-matrices are commutative, we have by rearranging (assuming $j \neq 2, p$)

$$\begin{aligned} v_i &\in \operatorname{col}\left[(A-\lambda_2 I)\cdots(A-\lambda_j I)\cdots(A-\lambda_p I)\right] \\ &= \operatorname{col}\left[(A-\lambda_j I)\cdot M\right], \text{ where } M = (A-\lambda_2 I)\cdots(A-\lambda_p I) \\ &= \operatorname{col}\left[(a_1\ a_2\ \cdots\ a_n)\cdot M\right] \\ &= \operatorname{col}\left(r_1 a_1 + r_2 a_2 + \cdots + r_n a_n\ \cdots\ q_1 a_1 + q_2 a_2 + \cdots + q_n a_n\right) \\ &= \operatorname{col}\left(a_1\ a_2\ \cdots\ a_n\right), \end{aligned}$$

where $a_1, a_2, \ldots, a_n$ are the columns of $(A-\lambda_j I)$, and $[r_1, \cdots, r_n]^T$ and $[q_1, \cdots, q_n]^T$ are the first and last columns of $M$, respectively. Other columns of $M$ are considered similarly. This, then implies that $v_i = k_1 a_1 + k_2 a_2 + \cdots + k_n a_n$ for some $k_1, k_2, \ldots, k_n \in \mathbb{R}$. Hence, $v_i$ is in the column space of $\kappa_{\lambda_j}(A)$. Similarly, we can show that $v_i$ is in all other $\kappa$-matrices for $j \neq i$. This completes the proof. □

According to Corollary 3.11, the complementary $\kappa$-matrices contain the eigenvectors as linear combinations of columns. That is, the eigenvector $v_i$ of $A$ for the eigenvalue $\lambda_i$ is contained in $\kappa_{\lambda_j}(A)$ for all $j \neq i$. $v_i$ may or may not be in $\operatorname{col}(\kappa_{\lambda_i}(A))$. It is also interesting to note that the eigenvalues of $A$ need not necessarily be in the spectrum of $\kappa_{\lambda_j}(A)$ for it to contain the associated eigenvectors.

The explanation above leads to the following result.

**Corollary 3.12** (Null-Column Theorem)**.** *Let $\lambda$ and $\mu$ be two eigenvalues of $A$, not necessarily distinct if they are repeated eigenvalues. Then* $\operatorname{Nul}(A-\lambda I) \subseteq \operatorname{Col}(A-\mu I)$.

Corollary 3.11 has a computation advantage. When the matrix has a spectrum of 2, we can take the full advantage of it. Following example is in order, which shows that eigenvectors already appear as columns of certain matrices.

**Example 3.13.** Let $\mathbf{A} = \begin{pmatrix} -1 & 1 & 1 \\ -4 & 3 & 2 \\ -6 & 3 & 4 \end{pmatrix}$, which has eigenvalues $\lambda = 4, 1$ and 1. The $\kappa$-matrix

$$\mathbf{A} - 1I = \begin{pmatrix} -2 & 1 & 1 \\ -4 & 2 & 2 \\ -6 & 3 & 3 \end{pmatrix}$$

has eigenvalues $\lambda = 3, 0, 0$ and contains, as nonzero column(s), the eigenvector of $A$ for $\lambda = 4$, which is not an eigenvalue of $(A - 1I)$. The $\kappa$-matrix

$$\mathbf{A} - 4I = \begin{pmatrix} -5 & 1 & 1 \\ -4 & -1 & 2 \\ -6 & 3 & 0 \end{pmatrix}$$

has eigenvalues $\lambda = 0, -3, -3$ and contains, as nonzero columns, the eigenvectors of $A$ for $\lambda = 1$, which is not an eigenvalue of $(A - 4I)$. The eigenvalues of $\kappa$-matrices can be obtained from the Shifting Theorem given in (3.1).

Let's consider several additional illustrative examples to fully understand the underlying argument.

3.3. **Special Methods for $3 \times 3$ Matrices.** There are several computational approaches when the diagonalizable matrix is of size 3. The underlying vector space can be taken as $\mathbb{R}$ for simplicity and hence, it can be related to the geometric representation of the 3-dimensional space. The case of spectrum of 1 is not interesting since any nonzero vector, $v_1$ can be made to an eigenvector by suitably selecting two other vectors to form an invertible matrix $P = (v_1\ v_2\ v_3)$ and such matrices are diagonal matrices with the eigenvalues sitting on the diagonal repeated three times to produce the equation $A = PDP^{-1}$, where $D = \operatorname{diag}(\lambda\ \lambda\ \lambda)$. In this case we assume that the eigenvalue is nonzero. Otherwise it generates a zero matrix as the product.

The case of spectrum of 2 has been considered earlier. So let's consider the case in which $A$ has three eigenvalues. Let them be $\lambda_1, \lambda_2$ and $\lambda_3$. According to The Eigenmatrix Theorem, the eigenvector corresponding to $\lambda_3$ lies in both $\kappa_{\lambda_1}(A)$ and $\kappa_{\lambda_2}(A)$. Let's look at a concrete example to fully understand the idea behind.

**Example 3.14.** Consider the $3 \times 3$ matrix $A$ with eigenvalues $\lambda_1 = 0, \lambda_2 = 3$ and $\lambda_3 = -4$ with the diagonalization given below.

$$\begin{pmatrix} 1 & 2 & 1 \\ 6 & -1 & 0 \\ -1 & -2 & -1 \end{pmatrix} = \begin{pmatrix} -1 & -2 & -1 \\ -6 & -3 & 2 \\ 13 & 2 & 1 \end{pmatrix} \begin{pmatrix} 0 & 0 & 0 \\ 0 & 3 & 0 \\ 0 & 0 & -4 \end{pmatrix} \begin{pmatrix} 1/12 & 0 & 1/12 \\ -8/21 & -1/7 & -2/21 \\ -9/28 & 2/7 & 3/28 \end{pmatrix} \tag{3.17}$$

The three $\kappa$-matrices are given by

$$\kappa_0(A) = \overset{\lambda_1 = 0}{\begin{bmatrix} 1 & 2 & 1 \\ 6 & -1 & 0 \\ -1 & -2 & -1 \\ u_1 & u_2 & u_3 \end{bmatrix}} \quad \kappa_3(A) = \overset{\lambda_2 = 3}{\begin{bmatrix} -2 & 2 & 1 \\ 6 & -4 & 0 \\ -1 & -2 & -4 \\ v_1 & v_2 & v_3 \end{bmatrix}} \quad \kappa_{-4}(A) = \overset{\lambda_3 = 3}{\begin{bmatrix} 5 & 2 & 1 \\ 6 & 3 & 0 \\ -1 & -2 & 3 \\ w_1 & w_2 & w_3 \end{bmatrix}} \tag{3.18}$$

There are several ways to proceed to find the associated eigenvectors. We discuss each method in detail below.

3.3.1. ***Method 1: Vector Calculus Approach I.*** First consider a vector calculus approach that uses the cross-product of two vectors. For further calculations we select vectors with least number of nonzero terms and ones that are easy to work with. Following that argument, we pick $u_2, u_3, v_2, v_3, w_2$ and $w_3$. Then we use the following argument. The eigenvector corresponds to $\lambda_3 = -4$ lies in the two planes generated by the column spaces of $\kappa_{\lambda_1=0}(A)$ and $\kappa_{\lambda_2=3}(A)$. Thus, it should be the line of intersection of those two planes. To find that, we first find normal vectors of the two planes and then find a line that is perpendicular to those two normal vectors. We may find all those normal lines using the cross-products in vector calculus.

A normal vector to the $\operatorname{col}(\kappa_{\lambda_1=0}(A))$ may be given by

$$\vec{n}_1 = v_2 \times v_3 = \begin{vmatrix} \mathbf{i} & \mathbf{j} & \mathbf{k} \\ -2 & 1 & 2 \\ 1 & 0 & -1 \end{vmatrix} = -\mathbf{i} + 0\mathbf{j} - \mathbf{k} = \begin{bmatrix} -1 \\ 0 \\ -1 \end{bmatrix}. \tag{3.19}$$

After rescaling we may take, $\vec{n}_1 = [1,\ 0,\ 1]^T$. Similarly we have $\vec{n}_2 = [8,\ 3,\ 2]^T$. These two normal vectors then produce a vector that is perpendicular to both of these vectors given by $\vec{l} = \vec{n}_1 \times \vec{n}_2 = [-3,\ 6,\ 3]^T = 3[-1,\ 2,\ 1]^T$. Thus an eigenvector of $\lambda_3$ may be given by $V_{\lambda_3} = [-1,\ 2,\ 1]^T$. Similarly we find other two eigenvectors.

A different but an interesting argument appears as an answer to the question 'Why can we find the eigenvector by cross multiplying two linearly independent row vectors in a matrix $A - \lambda I$ ?' among the discussions in [21] and has been taken almost verbatim to keep the original author's idea intact.

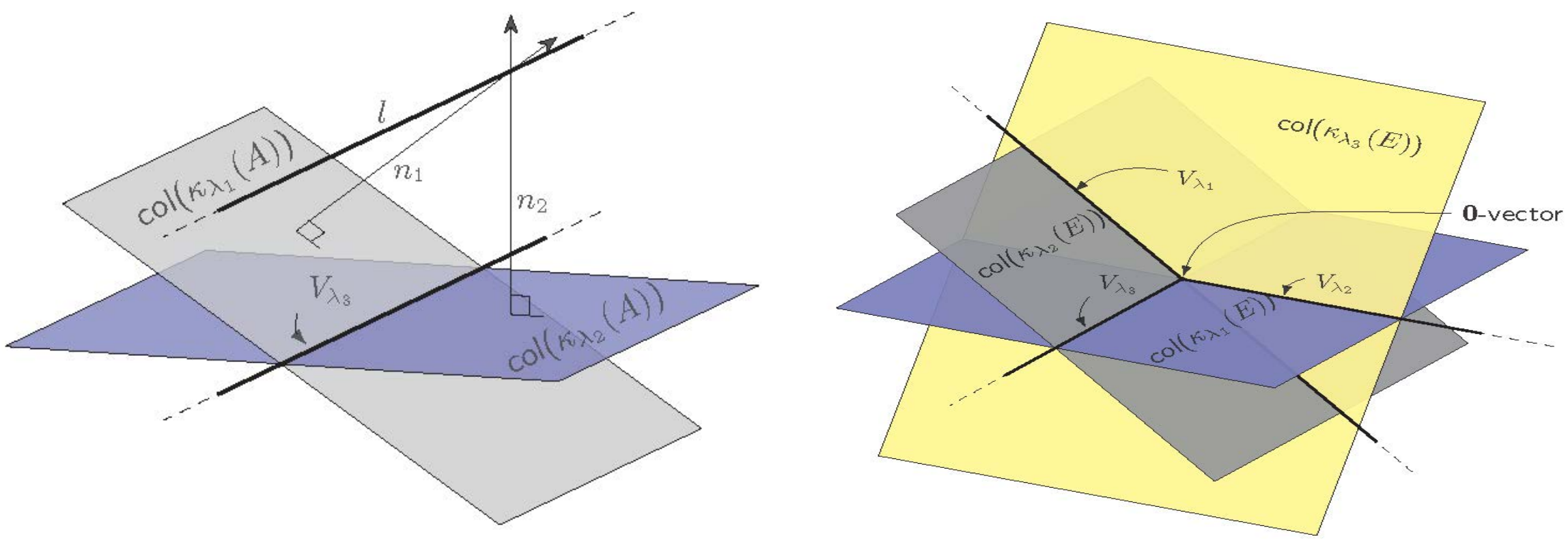

(A) Eigenspace of $\lambda_3$, $V_{\lambda_3}$ (B) Eigenspaces of $\lambda_1, \lambda_2$, and $\lambda_3$

FIGURE 5. Column spaces of $\kappa$-matrices

3.3.2. ***Method 2: Vector Calculus Approach II.*** The question is really about solving $M\mathbf{x} = \mathbf{0}$ when $M$ is a non-zero $3 \times 3$ matrix. Suppose $M$ has nonzero rows $\mathbf{M}_1, \mathbf{M}_2$, and $\mathbf{M}_3$, and suppose also that $\mathbf{M}_1 \nparallel \mathbf{M}_2$, so that $\mathbf{M}_1 \times \mathbf{M}_2 \neq \mathbf{0}$. Then the definition of the product of a matrix and a vector yields that

$$M(\mathbf{M}_1 \times \mathbf{M}_2) = \begin{bmatrix} \mathbf{M}_1 \cdot (\mathbf{M}_1 \times \mathbf{M}_2) \\ \mathbf{M}_2 \cdot (\mathbf{M}_1 \times \mathbf{M}_2) \\ \mathbf{M}_3 \cdot (\mathbf{M}_1 \times \mathbf{M}_2) \end{bmatrix} = \begin{bmatrix} 0 \\ 0 \\ \det M \end{bmatrix}.$$

If $M$ is singular, then $\det M = 0$ and so $\mathbf{M}_1 \times \mathbf{M}_2$ is a non-zero solution to $M\mathbf{x} = \mathbf{0}$. If two other non-parallel rows are used, then $M$ applied to their cross product will have two zero entries and $\pm \det M$ in the row corresponding to the vector triple product of all three rows of $M$, so we still get a non-zero solution to $M\mathbf{x} = \mathbf{0}$ when $M$ is singular.

Note that $M$ could be singular with all three rows being parallel, in which case $\ker M$ will be two-dimensional but cross products of row pairs will not give any non-zero results. Returning to the case of eigenvectors, we now have $M = \mathbf{A} - \lambda\mathbf{I}$ singular and so the cross product of two non-parallel rows of $\mathbf{A} - \lambda\mathbf{I}$ gives an eigenvector with eigenvalue $\lambda$. The case of all three rows being parallel will not occur if there are three distinct eigenvalues, since the eigenspace for each eigenvalue is one-dimensional and the observation just above requires a two-dimensional kernel.

3.3.3. ***Method 3: Matrix Product Approach.*** When $A$ is diagonalizable and has a spectrum of 2, each $\kappa$-matrix carries the eigenvectors of the complementary eigenvalue as nonzero columns. Following that argument, in the case of spectrum of 3, we combine last two $\kappa$-matrices to form one matrix, of which the columns are the eigenvectors of the complementary eigenvalue and consider other permutations of the matrices to produce eigenvectors of the other two eigenvalues. In the present example, the nonzero columns of the matrix obtained from the product

$$\kappa_3(A) \cdot \kappa_{-4}(A) = \begin{pmatrix} -2 & 2 & 1 \\ 6 & -4 & 0 \\ -1 & -2 & -4 \end{pmatrix} \begin{pmatrix} 5 & 2 & 1 \\ 6 & 3 & 0 \\ -1 & -2 & 3 \end{pmatrix} = \begin{pmatrix} 1 & 0 & 1 \\ 6 & 0 & 6 \\ -13 & 0 & -13 \end{pmatrix} \tag{3.20}$$

are the eigenvectors of $\lambda_1 = 0$. It is clear from (3.20) that the eigenvector is $[1,\ 6,\ -13]^T$ and it is easy to check that it is indeed an eigenvector of $A$ for the eigenvalue $\lambda_1 = 0$.

It is also worth to observe that we do not need to find the complete matrix product when finding the eigenvector(s) since we are only interested in one eigenvector. Thus, in practice, it is sufficient to multiply only one column of the matrix on the right with all the rows of the matrix on the left. If it is a zero vector, we need to select another column. The situation is illustrated below.

$$\kappa_3(A) \cdot \text{column 1 of } \kappa_{-4}(A) = \begin{pmatrix} -2 & 2 & 1 \\ 6 & -4 & 0 \\ -1 & -2 & -4 \end{pmatrix} \begin{pmatrix} 5 & * & * \\ 6 & * & * \\ -1 & * & * \end{pmatrix} = \begin{pmatrix} 1 & * & * \\ 6 & * & * \\ -13 & * & * \end{pmatrix}. \tag{3.21}$$

Since the resulting column is nonzero, it produces an eigenvector of $\lambda = 0$. Following this argument, we may find the eigenvector(s) of $\lambda_2 = 3$ by the product $\kappa_0(A) \cdot \kappa_{-4}(A)$ and of $\lambda_2 = -4$ by $\kappa_0(A) \cdot \kappa_3(A)$, respectively. Since there is only one eigenvector for each eigenvalue, we only need one such product to evaluate in each case. Thus, this approach will quickly find the necessary eigenvectors.

3.3.4. ***Method 4: Column Spaces Approach.*** According to the eigenmatrix theorem, each $\kappa$-matrix contains the eigenvectors of complimentary eigenvalues as linear combination of columns. Also, we only need to consider at most $n-1$ many nonzero columns in the process, due to the linear dependency of the columns of $\kappa$-matrices. We may use these facts to proceed as follows. Again in the present example, using two nonzero columns of $\kappa_3(A)$ and $\kappa_{-4}(A)$ each, we form a common linear combination as

$$a \begin{bmatrix} 2 \\ -4 \\ -2 \end{bmatrix} + b \begin{bmatrix} 1 \\ 0 \\ -4 \end{bmatrix} = c \begin{bmatrix} 2 \\ 3 \\ -2 \end{bmatrix} + d \begin{bmatrix} 1 \\ 0 \\ 3 \end{bmatrix} \tag{3.22}$$

The goal is to find $a, b, c$ and $d$ so that (3.22) produces a *common nonzero vector*. This then becomes an eigenvector of $\lambda_1 = 0$. We may use a trial-and-error approach or systematic approach such as the Gaussian-elimination at this step. For example, $a = 3, b = -8, c = -4$ and $d = 6$, gives $V_{\lambda=0} = [-2,\ -12,\ 26]^T$ as an eigenvector. Scaling the resulting vector, we may consider $v = [1,\ 6,\ -13]^T$ as an eigenvector. This method is sometimes faster than the Gaussian-elimination method.

3.3.5. ***Method 5: Trial-and-Error Approach I.*** Guess a common vector to $\kappa$-matrices of complementary eigenvalues. This common vector will then be an eigenvector of the complementary eigenvalue. This was the subject of one of the examples discussed earlier. The interested reader may refer back to Example 2.18 on page 10 for further details.

3.3.6. ***Method 6: Trial-and-Error Approach II.*** The $\kappa$-matrices most of the time provide potential candidates for eigenvectors as nonzero columns. Thus, we may use a trial-and-error approach to find the associated eigenvectors, possibly considering nonzero linear combinations as well. If all the choices fail, we may proceed with a other longer approaches.

3.3.7. ***Method 7: Gaussian-Elimination Approach.*** This is the classical approach to find the eigenvectors. The method works for any matrix in spite of being whether diagonalizable or not. This is the last resort if any of the above methods are not trivial or difficulty to work out. The method can be tedious when the order of the matrix becomes larger.

A second thought about why we never think of another way, or at least not known to the academic community, to find eigenvectors will be discussed in another section later.

According to the Diagonalization Theorem, a matrix with distinct eigenvalues, real or complex, is always diagonalizable. That means when a matrix is not diagonalizable, some of the eigenvalues should be repeated, possibly generating a case where it is difficulty to find an eigenbase for the underlying vector space. The next section will be devoted to study non-diagonalizable matrices in the context of eigenvectors.

## 4. Using Matrix Multiplication to Find Eigenvectors

As with the case of $3 \times 3$ matrices, we may use matrix multiplication to produce eigenvectors even for higher order matrices. However, in this case, it may be required to multiply several matrices before we may be able to generate any eigenvector. So, the time complexity may be an issue and the echelon matrix approach may be faster here. To compare the findings, we only consider one illustrative example here. Lets reconsider Example 2.22 that we solved using $\kappa$-matrix method.

**Example 4.1.** Again, let $H$ be given by

$$H = \begin{pmatrix} 3 & 2 & 5 & -5 \\ 3 & 4 & 7 & -7 \\ 4 & 4 & 10 & -9 \\ 6 & 6 & 14 & -13 \end{pmatrix}$$

As seen earlier, $H$ has three eigenvalues $\lambda_1 = 0$ and $\lambda_2 = 2$ each with multiplicity 1, and $\lambda_3 = 1$ with multiplicity 2. The corresponding eigenmatrices are given by

$$\kappa_0(H) = \underset{\lambda_1 = 0}{\begin{bmatrix} 3 & 2 & 5 & -5 \\ 3 & 4 & 7 & -7 \\ 4 & 4 & 10 & -9 \\ 6 & 6 & 14 & -13 \\ u_1 & u_2 & u_3 & u_4 \end{bmatrix}} \quad \kappa_2(H) = \underset{\lambda_2 = 2}{\begin{bmatrix} 1 & 2 & 5 & -5 \\ 3 & 2 & 7 & -7 \\ 4 & 4 & 8 & -9 \\ 6 & 6 & 14 & -15 \\ v_1 & v_2 & v_3 & v_4 \end{bmatrix}} \quad \kappa_1(H) = \underset{\lambda_3 = 1}{\begin{bmatrix} 2 & 2 & 5 & -5 \\ 3 & 3 & 7 & -7 \\ 4 & 4 & 9 & -9 \\ 6 & 6 & 14 & -14 \\ w_1 & w_2 & w_3 & w_4 \end{bmatrix}} \tag{4.1}$$

It can be seen that the product $\kappa_0(H) \cdot \kappa_2(H) \cdot \kappa_1(H) = 0$, which implies that $H$ is diagonalizable due to the Eigenmatrix Theorem 3.9. All we need to find is four linearly independent eigenvectors that form the transformation matrix $P$. In this example, we do that using the matrix product approach used for a $3 \times 3$ matrix earlier. According to the Method 3 of $3 \times 3$ matrices (See subsection 3.3.3), product of two $\kappa$-matrices will lead to an eigenvector of the complementary eigenvalue. For example, any nonzero column vector of the product $\kappa_2(H) \cdot \kappa_1(H)$ is an eigenvector associated with the eigenvalue $\lambda = 0$. Similar arguments work for other eigenvalues as well. This also happens as a result of the Eigenmatrix Theorem. However, we also know that we do not need to multiply the entire two matrices together. We only need to multiply the first matrix with any nonzero column vector of the second matrix. We normally use the one with least number of nonzero entries. The resulting product has to be either zero or an eigenvector of the complementary eigenvalue. For example, to obtain an eigenvector of $\lambda_1 = 0$, we multiply

$$\kappa_2(H) \cdot w_1 = \begin{pmatrix} -2 \\ -2 \\ -2 \\ -4 \end{pmatrix}$$

After rescaling, we take $\nu_1 = [1\ 1\ 1\ 2]^T$ as the eigenvector of $\lambda_1 = 0$. Similarly for $\lambda_2 = 2$, we use $\kappa_0(H) \cdot w_1 = [2\ 4\ 6\ 8]^T$, which gives the eigenvector $\nu_2 = [1\ 2\ 3\ 4]^T$. Finally for $\lambda_3 = 1$, we use $\kappa_0(H) \cdot v_1 = [-1\ 1\ 2\ 2]^T$, $\kappa_0(H) \cdot v_2 = [0\ 0\ 2\ 2]^T$, $\kappa_0(H) \cdot v_3 = [-1\ 1\ 2\ 2]^T$ and $\kappa_0(H) \cdot v_4 = [1\ -1\ -3\ -3]^T$, of which only two vectors are linearly independent. Considering two linear combination of the first two vectors, we may take $\nu_3 = [-1\ 1\ 0\ 0]^T$ and $\nu_4 = [0\ 0\ 1\ 1]^T$ as two linearly independent eigenvectors of $\lambda_3 = 1$, which has an algebraic multiplicity 2. It is worth noting that we may need to multiply each column vector of the second matrix, until it produces two linearly independent vectors. These two vectors then become the eigenvectors of the corresponding eigenvalue. The calculation also confirms that the geometric and algebraic multiplicities of each eigenvalue turns out to be equal, thus becomes a candidate for a diagonalizable matrix. Such a transformation can be given by

$$H = \begin{pmatrix} 1 & 1 & -1 & 0 \\ 1 & 2 & 1 & 0 \\ 1 & 3 & 0 & 1 \\ 2 & 4 & 0 & 1 \end{pmatrix} \begin{pmatrix} 0 & 0 & 0 & 0 \\ 0 & 2 & 0 & 0 \\ 0 & 0 & 1 & 0 \\ 0 & 0 & 0 & 1 \end{pmatrix} \begin{pmatrix} -1 & -1 & -3 & 3 \\ 1 & 1 & 2 & -2 \\ -1 & 0 & -1 & 1 \\ -2 & -2 & -2 & 3 \end{pmatrix}.$$

*Remark* 4.2. Now, there arises a complexity problem. According to the classical approach we first need to reduce the corresponding $\kappa$-matrix to the echelon form and then rewrite the corresponding system. This would then rearrange to produce the eigenvector(s). In the new approach, it only requires to multiply one column vector with one matrix. Thus, it seems that our approach is more efficient that the familiar Gaussian elimination method. In a future project, we will make an extended study to fully understand the time and computational complexities of the two approaches and will appear elsewhere.

## 5. What happens when a matrix is non-diagonalizable ?

The Eigenmatrix Theorem fails when a matrix is non-diagonalizable. However, the Equation 3.7 still holds even in this case if we consider all the eigenvalues taking into account the algebraic multiplicities and leads to an interesting result. To the contrary, the nonzero columns of product of $\kappa$-matrices produce generalized eigenvectors for eigenvalues with algebraic multiplicity more than 1. As a result, we will also consider a *generalized Jordan form* for such nondiagonalizable matrices in Section 7.3 later.

**Theorem 5.1** (Eigenmatrix Theorem for Non-diagonalizable Matrix)**.** *Let $\lambda_i, i = 1, \ldots, n$ be the eigenvalues, not necessarily distinct, of an $n \times n$ matrix $A$. Then*

$$(A - \lambda_1 I)(A - \lambda_2 I) \cdots (A - \lambda_n I) = 0 \tag{5.1}$$

*or*

$$\kappa_{\lambda_1}(A) \cdot \kappa_{\lambda_2}(A) \cdots \kappa_{\lambda_n}(A) = 0 \tag{5.2}$$

*In addition, the nonzero columns of the matrix*

$$\kappa_{\lambda_1}(A) \cdots \kappa_{\lambda_{i-1}}(A) \cdot \kappa_{\lambda_{i+1}}(A) \cdots \kappa_{\lambda_n}(A) \tag{5.3}$$

*are the eigenvectors of $A$ for the eigenvalue $\lambda_i, i = 2, \ldots, n-1$. The cases for $\lambda_i = 1$ and $n$ can be considered similarly.*

*Proof.* The proof is similar, if not equal, to the proof of the first part of Theorem 3.9 and reproduced here for completeness. Let $p(\lambda) = \lambda^n + a_{n-1}\lambda^{n-1} + \cdots + a_0$ be the characteristic polynomial of $A$. Then due to the Cayley-Hamilton Theorem for matrices, we have

$$p(A) = A^n + a_{n-1}A^{n-1} + \cdots + a_0 I = 0. \tag{5.4}$$

Suppose $A$ has $n$ eigenvalues $\lambda_i, i = 1, 2, \ldots, n$. Then by the Fundamental Theorem of Algebra, we may factor the left hand side of (5.4) as

$$(A - \lambda_1 I)(A - \lambda_2 I) \cdots (A - \lambda_n I) = 0. \tag{5.5}$$

Due to the commutativity of the $\kappa$-matrices, the order in which matrices are multiplied is not important. This completes the proof of the first part of the theorem. The second part follows a similar argument as in the second part of Theorem 3.9. This completes the proof. □

The eigenmatrix theorems for both diagonalizable and nondiagonalizable matrices are equally important and will produce eigenvectors when the totality of eigenvalues, taking into account the multiplicities, are considered since the product of $\kappa$-matrices being equal to zero means that the nonzero column vectors of the matrix on the right are the eigenvectors of the matrix produced as a result of the product of the complementary matrices in the product. This idea will be highly useful when the Jordan forms of nondiagonalizable matrices are concerned and will be discussed in the Section 7.3.

*Remark* 5.2. It is important to note in this case that the product $\kappa_{\lambda_1}(A)\cdot\kappa_{\lambda_2}(A)\cdots\kappa_{\lambda_p}(A) \neq 0$, if only the distinct eigenvalues are considered and leads to produce eigenvectors when multiplied by appropriate $\kappa$-matrices, not necessarily distinct. For example, if $\lambda_1$ has a multiplicity $m$ and $\lambda_2, \ldots, \lambda_p$ has multiplicities equal to 1, then the eigenvectors of $\lambda_1$ are the nonzero columns of the matrix produced by the product $\kappa_{\lambda_1}(A)\cdots\kappa_{\lambda_1}(A)\cdot\kappa_{\lambda_2}(A)\cdots\kappa_{\lambda_p}(A)$, where the first product has been taken $m-1$ many times. It can be seen that the number of nonzero and nonparallel, considered as a vector, columns is less than the algebraic multiplicity of the corresponding eigenvalue. That is, the geometric multiplicity is less than the algebraic multiplicity leading to a non-diagonalizable matrix.

In the next example, we find the eigenvectors of such a matrix using the argument in question. The example has been taken from [31].

**Example 5.3.** Let $\mathbf{A} = \begin{pmatrix} 2 & 4 & 3 \\ -4 & -6 & -3 \\ 3 & 3 & 1 \end{pmatrix}$. The eigenvalues are $\lambda_1 = 1$ and $\lambda_2 = -2$ with multiplicity 2. The corresponding $\kappa$-matrices are

$$\kappa_1 = \begin{pmatrix} 1 & 4 & 3 \\ -4 & -7 & -3 \\ 3 & 3 & 0 \end{pmatrix} \quad \text{and} \quad \kappa_{-2} = \begin{pmatrix} 4 & 4 & 3 \\ -4 & -4 & -3 \\ 3 & 3 & 3 \end{pmatrix} \tag{5.6}$$

It is easy to see that not all nonzero columns are eigenvectors of the complementary eigenvalue. However, we know that $\kappa_{\lambda_1}(A) \cdot \kappa_{\lambda_2}(A) \cdot \kappa_{\lambda_2}(A) = 0$, hence, as expected, the products $\kappa_{\lambda_2}(A) \cdot \kappa_{\lambda_2}(A)$ and $\kappa_{\lambda_1}(A) \cdot \kappa_{\lambda_2}(A)$ produce eigenvectors of $\lambda_1$ and $\lambda_2$, as nonzero columns, respectively. It can also be seen that the second product has not produced two eigenvectors, even though the associated eigenvalue has an algebraic multiplicity 2. As a result $A$ is not diagonalizable.

This gives rise to an additional question: what are those vectors produced as a product if they are not the eigenvectors of the complementary eigenvalues. It turns out these are the so-called *generalized eigenvectors* associated with the same set of eigenvalues. For example, we can see that $\kappa_{\lambda_2}(A) \cdot \kappa_{\lambda_2}(A) u_1 = \mathbf{0}$, despite the fact that $\kappa_{\lambda_2}(A) u_1 \neq \mathbf{0}$. This means that $u_1$ is a generalized eigenvector of $\lambda_2 = -2$, even though it is not an eigenvector of -2. To complete the discuss about $\kappa$-matrices, at least up to some extent, We will study generalized eigenvectors and possible *Jordan Canonical Forms* in a later section.

Thus, we may have a version of the diagonalization theorem as a criteria for testing whether a given matrix is diagonalizable. The corresponding shortcut formula is listed as the criteria 4 of the list of shortcut formulas in section 14. The result is an easy consequence of the Theorem 5.1.

## 6. Linear Transformations and $\kappa$-matrices

Linear transformations and matrices are closely related. Any linear transformation has a matrix representation with respect to a given basis. Let $V$ be a finite-dimensional vector space over $F$, where $F$ is either $\mathbb{R}$ or $\mathbb{C}$ depending on the context and $\mathcal{L}(V)$ be the set of linear operators on $V$. We start with the following definitions [5, p. 134–161].

**Definition 6.1** (Eigenvalue and Eigenvector)**.** Suppose $T \in \mathcal{L}(V)$. A scalar $\lambda \in F$ is called an *eigenvalue* of $T$ if there exists $\nu \in V$ such that $\nu \neq 0$ and $T\nu = \lambda\nu$. $\nu$ is called the *eigenvector* associated with $\lambda$.

In this work we only focus on finite-dimensional vector spaces $V$. The following result is standard in such vector spaces.

**Theorem 6.2.** *Every operator on a finite-dimensional, nontrivial, complex vector space has an eigenvalue.*

The core of the proof is based on the Fundamental Theorem of Algebra and factoring in $\mathbb{C}$ [5, p. 145]. Since we consider matrices in $\mathbb{R}^{n\times n}$ or $\mathbb{C}^{n\times n}$, we may equally apply the above result to obtain an eigenvalue for a given matrix $A$. Thus, the problems discussed throughout the paper are well-defined.

**Definition 6.3** (Matrix of an operator)**.** Let $T \in \mathcal{L}(V)$ and $\mathcal{B} = \{v_1, v_2, \ldots, v_n\}$ be a basis of $V$. The *matrix of* $T$ with respect to this basis is the $n \times n$ matrix

$$A = \begin{pmatrix} a_{11} & \cdots & a_{1n} \\ \vdots & & \vdots \\ a_{n1} & \cdots & a_{nn} \end{pmatrix}$$

whose entries $a_{jk}$ are given by

$$Tv_k = a_{1k}v_1 + \cdots + a_{nk}v_n.$$

Due to the above definition we can carry over all the results we had for matrices over $\mathbb{C}$ and do not distinguish unless otherwise stated. Not all $n \times n$ matrices have $n$ linearly independent eigenvectors. In such situations we move our attention to the so-called generalized eigenvectors. That is a subject of the next section. Further studies of the connection between linear operators and their $\kappa$-matrices will be studied in another project.

## 7. Generalized Eigenvectors of a Matrix

Some matrices cannot be diagonalized. It mainly depends on whether the matrix in question has $n$ linearly independent eigenvectors to form an eigenbasis for $\mathbb{R}^n$. In this section, we study the question why some matrices do not have enough eigenvectors. We also have a look at the more generalized Jordan Canonical Forms (JCF). A more detailed exposition about JCF will be discussed in a future project.

**Definition 7.1** (Generalized eigenvector [5, p. 245])**.** Suppose $T \in \mathcal{L}(V)$ and $\lambda$ is an eigenvalue of $T$. A vector $v \in V$ is called an *generalized eigenvector of rank* $j$ of $T$ corresponding to $\lambda$ if $v \neq \mathbf{0}$ and

$$(T - \lambda I)^j \, v = 0 \quad \text{and} \quad (T - \lambda I)^{j-1} \, v \neq 0$$

for some $j \in \mathbb{N}$.

According to this definition, the eigenvectors become rank 1 generalized eigenvectors of the operator $T$. The same is valid for matrices as well. It can be shown that $j \leq \dim V$. We have the following results.

**Theorem 7.2** (Linearly independent generalized eigenvectors [5, p. 247])**.** *Let $T \in \mathcal{L}(V)$. Suppose $\lambda_1, \ldots, \lambda_m$ are distinct eigenvalues of $T$ and $v_1, \ldots, v_m$ are corresponding generalized eigenvectors. Then $v_1, \ldots, v_m$ is linearly independent.*

*Remark* 7.3. When $T$ is an $n \times n$ matrix, it may have $m <= n$ eigenvalues with $n$ generalized eigenvectors. Thus, $T$ can be written in the JCF.

**Definition 7.4** (Generalized eigenspaces)**.** Let $T \in \mathcal{L}(V)$ and $\lambda \in F$. The *generalized eigenspace* of $T$ corresponding to $\lambda$, denoted by $G(\lambda, T)$, is defined to be the set of all generalized eigenvectors of $T$ corresponding to $\lambda$, along with the 0 vector.

Thus, the model matrix $P$ formed by the generalized eigenvectors is always invertible, leading to the Jordan Form of the matrix decomposition.

In other words, it contains all vectors, $\nu$ such that $(T - \lambda I)^j \nu = 0$ for some integer $j \geq 1$. However, it turns out that $j \leq m$, where $m$ is the algebraic multiplicity of $\lambda$. It can be shown that there exist $m$ for such linearly independent generalized eigenvectors. We have the following result for generalized eigenvectors [5] .

**Theorem 7.5.** *Suppose $V$ is a complex vector space and $T \in \mathcal{L}(V)$. Then there is a basis of $V$ consisting of generalized eigenvectors of $T$.*

This follows from the decomposition theorem of a vector space in terms of generalized eigenspaces and choosing a basis for each generalized eigenspace, which has a dimension same as the algebraic multiplicity of the associated eigenvalue.

For clarity and further development we may restate the following version of the generalized vector for a matrix $A$.

**Definition 7.6** (Generalized eigenvector [9, p. 189])**.** A vector $x_m$ is a *generalized eigenvector of rank* $m$ corresponding to the matrix $A$ and the eigenvalue $\lambda$ if $(A - \lambda I)^m \, x_m = 0$ but $(A - \lambda I)^{m-1} \, x_m \neq 0$.

We will use this definition for further discussions.

7.1. **Relation between generalized eigenvectors and $\kappa$-matrices.** We consider an example to discuss the argument behind. It has been taken from [31, p. 286].

**Example 7.7.** Let $A = \begin{pmatrix} 2 & 4 & 3 \\ -4 & -6 & -3 \\ 3 & 3 & 1 \end{pmatrix}$

$A$ has the eigenvalues $\lambda_1 = 1$ with multiplicity 1 and $\lambda_2 = -2$ with multiplicity 2. It can be seen that the product of the $\kappa$-matrices is

$$\kappa_1(A) \cdot \kappa_{-2}(A) = \begin{pmatrix} 1 & 4 & 3 \\ -4 & -7 & -3 \\ 3 & 3 & 0 \end{pmatrix} \begin{pmatrix} 4 & 4 & 3 \\ -4 & -4 & -3 \\ 3 & 3 & 3 \end{pmatrix} = \begin{pmatrix} -3 & -3 & 0 \\ 3 & 3 & 0 \\ 0 & 0 & 0 \end{pmatrix}$$

Thus, according to the 2-Spectrum lemma, $A$ is not diagonalizable. In addition, the $3^{rd}$ column of $\kappa_{-2}(A)$, that is $v_1 = [3, \ -3, \ 3]^T$ is an eigenvector of $\lambda_1 = 1$ and similarly considering $\kappa_{-2}(A) \cdot \kappa_1(A)$, $v_2 = [3, \ -3, \ 0]^T$ is an eigenvector of $\lambda_2 = -2$. It can be seen that there are no other linearly independent eigenvectors, thus making it difficult to bring it to the diagonal form. However, it has a Jordan form and we will discuss the theory in detail.

We need the definition below [9].

**Definition 7.8** (Chains)**.** Let $x_m$ be a generalized eigenvector of rank $m$ corresponding to the matrix $A$ and the eigenvalue $\lambda$. The *chain generated by* $x_m$ is a set of vectors $\{x_m, \ x_{m-1}, \ \ldots, \ x_1\}$ given by

$$\begin{aligned} x_{m-1} &= (\mathbf{A} - \lambda \mathbf{I})\, x_m \\ x_{m-2} &= (\mathbf{A} - \lambda \mathbf{I})^2\, x_m = (\mathbf{A} - \lambda \mathbf{I})\, x_{m-1} \\ x_{m-3} &= (\mathbf{A} - \lambda \mathbf{I})^3\, x_m = (\mathbf{A} - \lambda \mathbf{I})\, x_{m-2} \\ &\vdots \\ x_1 &= (\mathbf{A} - \lambda \mathbf{I})^{m-1}\, x_m = (\mathbf{A} - \lambda \mathbf{I})\, x_2. \end{aligned}$$

In general,

$$x_j = (\mathbf{A} - \lambda \mathbf{I})^{m-j}\, x_m = (\mathbf{A} - \lambda \mathbf{I})\, x_{j+1}, (j = 1, 2, \ldots, m-1). \tag{7.1}$$

The prove of the following important theorem can be found, for example in [9, p. 194–196].

**Theorem 7.9.** *$x_j$ given by* (7.1) *is a generalized vector of rank $j$ corresponding to the eigenvalue $\lambda$ and a chain $\{x_m, \ x_{m-1}, \ \ldots, \ x_1\}$ is a linearly independent set of vectors.*

**Definition 7.10** (Canonical basis)**.** A set of $n$ linearly independent generalized eigenvectors is a *canonical basis* if it is composed of entirely of chains.

Now, we are ready to look at the Jordan Canonical Forms for both diagonalizable and nondiagonalizable matrices.

7.2. **Jordan Canonical Forms.** We refer to the classical books, for example, [9, 12, 18, 36] for further details.

According to the classical theory of Jordan canonical forms, every matrix is similar to an "almost diagonal" matrix. For rest of the definitions of this section we follow [9, Section 8.7].

Let $\lambda_k$ be an eigenvalue of $A$. Now define a matrix $S_k$ that has all of its diagonal elements equal to $\lambda_k$, all of its superdiagonal elements equal to 1, and the rest of the elements equal to zero, given by

$$\mathbf{S}_k = \begin{pmatrix} \lambda_k & 1 & 0 & 0 & \cdots & 0 & 0 \\ 0 & \lambda_k & 1 & 0 & \cdots & 0 & 0 \\ 0 & 0 & \lambda_k & 1 & \cdots & 0 & 0 \\ \vdots & \vdots & \vdots & \vdots & \ddots & \vdots & \vdots \\ 0 & 0 & 0 & 0 & \cdots & \lambda_k & 1 \\ 0 & 0 & 0 & 0 & \cdots & 0 & \lambda_k \end{pmatrix} \tag{7.2}$$

**Definition 7.11** (Jordan Canonical Form)**.** [9, p. 203–204] A square matrix $A$ is in *Jordan canonical form* if it is a diagonal matrix or can be written in either one of the following forms:

$$\begin{pmatrix} \mathbf{D} & & & & \\ & \mathbf{S}_1 & & & \\ & & \cdot & & \\ & & & \cdot & \\ & & & & \cdot \\ & & & & & \mathbf{S}_r \end{pmatrix}$$

or

$$\begin{pmatrix} \mathbf{S}_1 & & & & \\ & \cdot & & & \\ & & \cdot & & \\ & & & \cdot & \\ & & & & \mathbf{S}_r \end{pmatrix}$$

Here $\mathbf{D}$ is a diagonal matrix and $\mathbf{S}_k (k = 1, 2, \ldots, r)$ is defined by (7.2).

It turns out that when the "generalized eigenvectors" are not part of the chains, they do not satisfy the classical JCF, rather it satisfies a more general form. We focus on to that topic in the next section.

7.3. **Variation of the Jordan Canonical Forms.** We start with the following observation.

**Observation 7.12.** Notice that in Example 7.7, we have the *Jordan-like* decomposition given by

$$A = \begin{pmatrix} 3 & 3 & 1 \\ -3 & -3 & -4 \\ 3 & 0 & 3 \end{pmatrix} \begin{pmatrix} 1 & 0 & 0 \\ 0 & -2 & -1 \\ 0 & 0 & -2 \end{pmatrix} \begin{pmatrix} 1/3 & 1/3 & 1/3 \\ 4/9 & 1/9 & -1/3 \\ -1/3 & -1/3 & 0 \end{pmatrix} \tag{7.3}$$

where the model matrix $P$ consists of the $3^{rd}$ column of $\kappa_{-2}(A)$, $3^{rd}$ column of $\kappa_1(A)$, and $2^{nd}$ column of $\kappa_1(A)$, as eigenvectors/generalized eigenvectors of $A$.

Now, on the other hand, if we switch the $2^{nd}$ and $3^{rd}$ columns of $P$ we obtain

$$A = \begin{pmatrix} 3 & 1 & 3 \\ -3 & -4 & -3 \\ 3 & 3 & 0 \end{pmatrix} \begin{pmatrix} 1 & 0 & 0 \\ 0 & -2 & 0 \\ 0 & -1 & -2 \end{pmatrix} \begin{pmatrix} 1/3 & 1/3 & 1/3 \\ -1/3 & -1/3 & 0 \\ 4/9 & 1/9 & -1/3 \end{pmatrix} \tag{7.4}$$

Note that the diagonal-like matrix is the transpose of the diagonal-like matrix of the previous case in Eq. 7.3. This may not come as a surprise since the rearrangement of the vectors in the model matrix (transformation matrix) requires some short of rearrangement of entries in the corresponding diagonal matrix as well.

It can be seen that we may use the $1^{st}$ column of $\kappa_1(A)$ to arrive at the same answer as well.

**Observation 7.13.** Consider the matrix $A$, which has three distinct eigenvalues $\lambda_1 = 1, \lambda_2 = 2$ and $\lambda_3 = 4$ with multiplicity 2. It can be seen that $A$ is not diagonalizable due to the eigenmatrix theorem.

$$A = \begin{pmatrix} 5 & 4 & 2 & 1 \\ 0 & 1 & -1 & -1 \\ -1 & -1 & 3 & 0 \\ 1 & 1 & -1 & 2 \end{pmatrix}$$

which has the $\kappa$-matrices

$$\kappa_1(A) = \begin{pmatrix} 4 & 4 & 2 & 1 \\ 0 & 0 & -1 & -1 \\ -1 & -1 & 2 & 0 \\ 1 & 1 & -1 & 1 \end{pmatrix}, \quad \kappa_2(A) = \begin{pmatrix} 3 & 4 & 2 & 1 \\ 0 & -1 & -1 & -1 \\ -1 & -1 & 1 & 0 \\ 1 & 1 & -1 & 0 \end{pmatrix}, \quad \kappa_4(A) = \begin{pmatrix} 1 & 4 & 2 & 1 \\ 0 & -3 & -1 & -1 \\ -1 & -1 & -1 & 0 \\ 1 & 1 & -1 & -2 \end{pmatrix}$$

Considering the products of complementary $\kappa$-matrices and rescaling, we produce $v_{\lambda=1} = [1,\ -1,\ 0,\ 0]^T$, $v_{\lambda=2} = [1,\ -1,\ 0,\ 1]^T$ and $v_{\lambda=4} = [1,\ 0,\ -1,\ 1]^T$ as eigenvectors of 1, 2 and 4, respectively. To find the generalized eigenvector of 4 of rank 2, we consider the product $\kappa_1(A) \cdot \kappa_2(A)$ and obtain

$$\kappa_1(A) \cdot \kappa_2(A) = \begin{pmatrix} 11 & 11 & 5 & 0 \\ 0 & 0 & 0 & 0 \\ -5 & -5 & 1 & 0 \\ 5 & 5 & -1 & 0 \end{pmatrix}$$

Now there are two candidates for the generalized vector and it can be shown that they are indeed rank 2 generalized eigenvectors of 4. That is, taking $v_4 = [11,\ 0,\ -5,\ 5]^T$ and $v_5 = [5,\ 0,\ 1,\ -1]^T$, we see that $\kappa_4(A) \cdot \kappa_4(A) v_4 = \kappa_4(A) \cdot \kappa_4(A) v_5 = 0$. Let's use them to produce Jordan-like forms. First taking $v_4 = [11,\ 0,\ -5,\ 5]^T$ as the generalized eigenvector of 4, we have

$$A = [v_{\lambda=1}\ v_{\lambda=2}\ v_{\lambda=4}\ v_4] \begin{pmatrix} 1 & 0 & 0 & 0 \\ 0 & 2 & 0 & 0 \\ 0 & 0 & 4 & \boxed{6} \\ 0 & 0 & 0 & 4 \end{pmatrix} [v_{\lambda=1}\ v_{\lambda=2}\ v_{\lambda=4}\ v_4]^{-1}. \tag{7.5}$$

It can be seen that the choice of $v_5$ as the generalized vector produce the same result.

Now if we use the nonrescaled eigenvector of $\lambda = 4$, that is $v_{\lambda=4} = [6,\ 0,\ -6,\ 6]^T$ we have

$$A = [v_{\lambda=1}\ v_{\lambda=2}\ v_{\lambda=4}\ v_4] \begin{pmatrix} 1 & 0 & 0 & 0 \\ 0 & 2 & 0 & 0 \\ 0 & 0 & 4 & \boxed{1} \\ 0 & 0 & 0 & 4 \end{pmatrix} [v_{\lambda=1}\ v_{\lambda=2}\ v_{\lambda=4}\ v_4]^{-1}. \tag{7.6}$$

*Remark* 7.14. Thus, it can be seen that *we cannot rescaled* generalized eigenvectors of associated with a particular eigenvalue. If we would like to rescale, we need to rescale them all within one chain. Now if we use the classical definition of the generalized vector, that is Equation 7.1, we have

$$v_{\lambda=4} = (A - 4I)\, v_{gv}$$

This would then give $v_{gv} = [6\ 0\ 0\ 0]^T$ or $v_{gv} = [0\ 0\ 6\ -6]^T$ as the rank 2 generalized vector. However, these vectors cannot be scaled if we are to use them in the Jordan chains.

As the next example, we consider the $5 \times 5$ matrix below that appears in [20].

**Observation 7.15.**

$$A = \begin{pmatrix} 1 & 0 & 0 & 0 & 0 \\ 3 & 1 & 0 & 0 & 0 \\ 6 & 3 & 2 & 0 & 0 \\ 10 & 6 & 3 & 2 & 0 \\ 15 & 10 & 6 & 3 & 2 \end{pmatrix}$$

which has eigenvalues $\lambda_1 = 1$ and $\lambda_2 = 2$ with algebraic multiplicities 2 and 3 and geometric multiplicities 1.

We do not follow the classical chain approach to find the generalized vectors, but rather a matrix-vector product approach. So we first consider the corresponding $\kappa$-matrices:

$$\kappa_1(A) = \begin{pmatrix} 0 & 0 & 0 & 0 & 0 \\ 3 & 0 & 0 & 0 & 0 \\ 6 & 3 & 1 & 0 & 0 \\ 10 & 6 & 3 & 1 & 0 \\ 15 & 10 & 6 & 3 & 1 \end{pmatrix}, \quad \kappa_2(A) = \begin{pmatrix} -1 & 0 & 0 & 0 & 0 \\ 3 & -1 & 0 & 0 & 0 \\ 6 & 3 & 0 & 0 & 0 \\ 10 & 6 & 3 & 0 & 0 \\ 15 & 10 & 6 & 3 & 0 \end{pmatrix}.$$

First, it can be seen that $\kappa_1(A)\cdot\kappa_2(A)\neq 0$ meaning that $A$ is not diagonalizable. So, to find the generalized eigenvectors we use the following approach.

Consider the following matrix products:

$$\kappa_2(A)\cdot\kappa_2(A)\cdot\kappa_2(A)\cdot\kappa_1(A)=\begin{pmatrix}0&0&0&0&0\\-3&0&0&0&0\\9&0&0&0&0\\-9&0&0&0&0\\3&0&0&0&0\end{pmatrix},$$
$$\kappa_2(A)\cdot\kappa_2(A)\cdot\kappa_2(A)=\begin{pmatrix}-1&0&0&0&0\\9&-1&0&0&0\\-12&3&0&0&0\\-17&-3&0&0&0\\51&1&0&0&0\end{pmatrix} \tag{7.7}$$

The *first columns* of these two matrices are the generalized vectors of $\lambda_1=1$, which has a multiplicity 2. Now, for the other eigenvalue consider the following three matrix products.

$$\kappa_1(A)\cdot\kappa_1(A)\cdot\kappa_2(A)\cdot\kappa_2(A)=\begin{pmatrix}0&0&0&0&0\\0&0&0&0&0\\0&0&0&0&0\\0&0&0&0&0\\135&27&9&0&0\end{pmatrix},$$
$$\kappa_1(A)\cdot\kappa_1(A)\cdot\kappa_2(A)=\begin{pmatrix}0&0&0&0&0\\0&0&0&0&0\\0&0&0&0&0\\45&9&3&0&0\\228&63&24&3&0\end{pmatrix}, \tag{7.8}$$
$$\kappa_1(A)\cdot\kappa_1(A)=\begin{pmatrix}0&0&0&0&0\\0&0&0&0&0\\15&3&1&0&0\\46&15&6&1&0\\111&46&21&6&1\end{pmatrix}.$$

The *first columns* of these three matrices are the generalized vectors of $\lambda_2=2$, which has a multiplicity 3. It is an interesting observation that the *second* or *third* columns of these matrices are also the generalized eigenvectors of 2 in a single chain, as long as we only consider the same columns of each matrix. The reason for that will be studied in detail in a future project.

Now taking the generalized vectors in question as the vectors for the transformation matrix we may produce the Jordan form as follows, even though these vectors are not necessarily the vectors discussed in Jordan chains in the classical sense.

$$A=\begin{pmatrix}0&-1&0&0&0\\-3&9&0&0&0\\9&-12&0&0&15\\-9&-17&0&45&46\\3&51&135&228&111\end{pmatrix}\begin{pmatrix}1&1&0&0&0\\0&1&0&0&0\\0&0&2&1&0\\0&0&0&2&1\\0&0&0&0&2\end{pmatrix}(\cdot)^{-1}.$$

where $(\cdot)$ is the first matrix appeared in the product. As mentioned earlier, we may scaled the vectors within a chain of one eigenvalue. Also, in practice we tend to use a chain with nonzero vectors that has many zeros and that are easy to work with. For example, in this example, it would have been easier if we had adopted the third columns for the eigenvalue $\lambda_2=2$, instead, may be also scaled the chain by 3, that means dividing each element in the chain by 3.

As the final example regarding the Jordan forms, we consider the following $5\times 5$ matrix discussed in [36, p. 123].

**Observation 7.16.** Let

$$A = \begin{pmatrix} 1 & 0 & -1 & 1 & 0 \\ -4 & 1 & -3 & 2 & 1 \\ -2 & -1 & 0 & 1 & 1 \\ -3 & -1 & -3 & 4 & 1 \\ -8 & -2 & -7 & 5 & 4 \end{pmatrix}$$

which has $\lambda = 2$ as the eigenvalue with algebraic multiplicity 5 and geometric multiplicity 2 since the Hermite normal form (or the reduced-echelon form) of $(A - 2I)$ has two free variables corresponding to the last two columns of $(A - 2I)$. So we consider the chains associated with these two columns, one of which has a length 2 and the other has length 3. To do that we first consider the products (powers) of $\kappa$-matrices below:

$$\kappa_2(A) = \begin{pmatrix} -1 & 0 & -1 & 1 & 0 \\ -4 & -1 & -3 & 2 & 1 \\ -2 & -1 & -2 & 1 & 1 \\ -3 & -1 & -3 & 2 & 1 \\ -8 & -2 & -7 & 5 & 2 \end{pmatrix}, \tag{7.9}$$

$$(\kappa_2(A))^2 = \begin{pmatrix} 0 & 0 & 0 & 0 & 0 \\ 0 & 0 & 0 & 0 & 0 \\ -1 & 0 & -1 & 1 & 0 \\ -1 & 0 & -1 & 1 & 0 \\ -1 & 0 & -1 & 1 & 0 \end{pmatrix}, \tag{7.10}$$

and $(\kappa_2(A))^3 = 0$. Thus the nonzero columns of $(A - 2I)$ are the generalized eigenvectors of rank 2.

To find the second element of the chain of length 2, take $\nu_1 = (0,\ 1,\ 1,\ 1,\ 2)$ and solve

$$(A - 2I)\, x_2 = \nu_1$$

for $x_2$, which is the second element of the chain. This then gives $x_2 = (0,\ -1,\ 0,\ 0,\ 0) + t_1\,(0,\ -1,\ 1,\ 1,\ 0) + t_2\,(0,\ 1,\ 0,\ 0,\ 1)\,,\ t_1, t_2 \in \mathbb{R}$. Taking $t_1 = 0$ and $t_2 = 1$, we have, in particular, $x_2 = (0,\ 0,\ 0,\ 0,\ 1)$.

Similarly, to find the third element, $x_3$ of the chain of length 3, let $\nu_2 = (1,\ 2,\ 1,\ 2,\ 5)$ and obtain $x_3 = (0,\ 1,\ -1,\ 0,\ 0) + t_1\,(0,\ -1,\ 1,\ 1,\ 0) + t_2\,(0,\ 1,\ 0,\ 0,\ 1)\,,\ t_1, t_2 \in \mathbb{R}$. Taking $t_1 = 1$ and $t_2 = 0$, we have arrive at $x_3 = (0,\ 0,\ 0,\ 1,\ 0)$. The other choices would work just as well.

This then leads to the Jordan form

$$A = PDP^{-1} = \begin{pmatrix} 0 & 0 & 0 & 1 & 0 \\ 1 & 0 & 0 & 2 & 0 \\ 1 & 0 & 1 & 1 & 0 \\ 1 & 0 & 1 & 2 & 1 \\ 2 & 1 & 1 & 5 & 0 \end{pmatrix} \begin{pmatrix} 2 & 1 & 0 & 0 & 0 \\ 0 & 2 & 0 & 0 & 0 \\ 0 & 0 & 2 & 1 & 0 \\ 0 & 0 & 0 & 2 & 1 \\ 0 & 0 & 0 & 0 & 2 \end{pmatrix} (\cdot)^{-1}.$$

This method of obtaining the Jordan forms is faster than the classical approach. The generalized vectors that form the chains are simply the vectors in the same columns of the sequences of $\kappa$-matrix products. The final element in the chain can be obtain by following the classical approach as was the case in the example in question.

Leaving room for further development, the extended studies on the generalized Jordan forms will be done in a future project.

Next we answer the following philosophical question.

## 8. Why we did it the other way in the first place?

According to the natural flow of thinking, we tend to use the Gaussian elimination (echelon forms) when it comes to solve a system of linear equations. Finding eigenvectors boil down to finding the basis elements of nullspaces, that is solving the matrix equation $(A - \lambda I)x = 0$. The most natural and obvious approach was to use the Gaussian elimination, thus leaving no need for another method. The argument goes as follows.

Consider the problem of finding the eigenvectors of the matrix $A$ for the eigenvalue $\lambda$. Then we use the following method to find them.

- **Procedure: The echelon form approach to find eigenvectors**
  1) Write the augmented matrix $\left[A-\lambda I \middle| 0\right]$.
  2) Use *Gaussian elimination* to row-reduce the matrix to reduced echelon form.
  3) Write the corresponding system of equation using nonzero rows.
  4) Solve for basic variables in terms of free variables.
  5) Find the basis elements of $\mathrm{Nul}(A-\lambda I)$

The basis element(s) then become the eigenvector(s) associated with the eigenvalue $\lambda$.

The approach is simple and straightforward. The implementation on a computer algebra system (CAS) is again easy and can also use systematic elementary matrix approach at times. CAS such as Octave®, Matlab®, Mathematica® and Maple® have built-in subroutines to implement the classical Gaussian approach with some optimization at hand.

As the final illustration example, we consider the following matrix, which has both complex and non-integer entries to make it extremely difficulty solve using a paper-pencil approach.

**Example 8.1.** Consider the following $5\times 5$ real matrix with three distinct real and two complex eigenvalues. Let $K$ be given by

$$K=\begin{pmatrix} 1-1.5i & 2-i & -4-2.5i & 1+3i & -3 \\ -2+i & -4+i & 3 & 1-i & 5-i \\ -1+1.5i & i & 2+2.5i & -1-3i & 1 \\ -2+i & -1+i & -1 & 2-i & 1-i \\ -1-0.5i & -3 & -2.5i & 2+2i & 3-i \end{pmatrix}$$

It has five eigenvalues each with multiplicity 1, namely, $\lambda_1=0, \lambda_2=1, \lambda_3=-1, \lambda_4=2-i$ and $\lambda_5=2+i$. According to the argument used in Example 2.22, we only need to consider two $\kappa$-matrices for two different eigenvalues. For simplicity, we use $\lambda_1=0$ and $\lambda_2=1$. The corresponding eigenmatrices are given by

$$\overset{\kappa_{\lambda_1}(K)}{\begin{bmatrix} 1-1.5i & 2-i & -4-2.5i & 1+3i & -3 \\ -2+i & -4+i & 3 & 1-i & 5-i \\ -1+1.5i & i & 2+2.5i & -1-3i & 1 \\ -2+i & -1+i & -1 & 2-i & 1-i \\ -1-0.5i & -3 & -2.5i & 2+2i & 3-i \\ u_1 & u_2 & u_3 & u_4 & u_5 \end{bmatrix}}\overset{\kappa_{\lambda_2}(K)}{\begin{bmatrix} -1.5i & 2-i & -4-2.5i & 1+3i & -3 \\ -2+i & -5+i & 3 & 1-i & 5-i \\ -1+1.5i & i & 1+2.5i & -1-3i & 1 \\ -2+i & -1+i & -1 & 1-i & 1-i \\ -1-0.5i & -3 & -2.5i & 2+2i & 2-i \\ v_1 & v_2 & v_3 & v_4 & v_5 \end{bmatrix}}$$

We then use $u_2,\ldots,u_5, v_2,\ldots,v_5$ for further calculations. The eigenvectors of complimentary eigenvalues take the form

$$\begin{aligned} v=au_2+bu_3+cu_4+du_5 &= a\begin{pmatrix} 2-i \\ -4+i \\ i \\ -1+i \\ -3 \end{pmatrix} + b\begin{pmatrix} -4-2.5i \\ 3 \\ 2+2.5i \\ -1 \\ -2.5i \end{pmatrix} + c\begin{pmatrix} 1+3i \\ 1-i \\ -1-3i \\ 2-i \\ 2+2i \end{pmatrix} + d\begin{pmatrix} -3 \\ 5-i \\ 1 \\ 1-i \\ 3-i \end{pmatrix} \\ &= \begin{pmatrix} 2a-4b+c-3d+(-a-2.5b+3c)i \\ -4a+3b+c+5d+(a-c-d)i \\ 2b-c+d+(a+2.5b-3c)i \\ -a-b+2c+d+(a-c-d)i \\ -3a+2c+3d+(-2.5b+2c-d)i \end{pmatrix}. \end{aligned} \tag{8.1}$$

Now, for $\lambda_2=1$, we have from $Kv=v$ and the fourth row of $K$ that

$$\begin{aligned} (-2+i)\left[2a-4b+c-3d+(-a-2.5b+3c)i\right] + (-1+i)\left[-4a+3b+c+5d+(a-c-d)i\right] \\ +(-1)\left[2b-c+d+(a+2.5b-3c)i\right] + (2-i)\left[-a-b+2c+d+(a-c-d)i\right] \\ +(1-i)\left[-3a+2c+3d+(-2.5b+2c-d)i\right] = -a-b+2c+d+(a-c-d)i \end{aligned}$$

This gives

$$-3a + 2b + c + 3d + 3(a - c - d)i = 0 \Longrightarrow a = c + d \text{ and } b = c, \text{where } c, d \in \mathbb{R}. \tag{8.2}$$

This implies that $c$ and $d$ can be any real numbers. Now replacing $a$ with $c + d$ and $b$ with $c$, we have by (8.1) that

$$v = \begin{pmatrix} -c - d + (-0.5c - d)i \\ d \\ c + d + (0.5c + d)i \\ 0 \\ -c + (-0.5c - d)i \end{pmatrix} \tag{8.3}$$

Now, using the fifth rows of $K$ and $v$, we have

$$(-1 - 0.5i)(-c - d + (-0.5c - d)i) - 3d - 2.5i[c + d + (0.5c + d)i] + 0 + (3 - i)[-c + (-0.5c - d)i] \\ = -c + (-0.5c - d)i,$$

which implies that $c = -2d$, where $d \in \mathbb{R}$. Since the nullspace of $\lambda = 1$ has dimension 1, it is spanned by a vector $v$ given by $v = (1,\ 1,\ -1,\ 0,\ 2)^T$, which is an eigenvector of 1. Also note that we only need to consider two rows of the eigenformula $Kv = v$ since, by doing so we can express $v$ as an one parameter family.

The Gauss-Jordan elimination of $(K - 1 \cdot I)$ can be lengthy and tedious due to the involvement of non-integer entries and complex numbers. However, the present method provides a venue for easy path to finding such complicated eigenvectors with a relatively less complex approach. The other eigenvectors can be found similarly and left as an exercise for interested readers.

We conclude the major results of our work with the following observation that comes as a result of left-multiplication of row vectors and matrices.

## 9. Right vs Left Eigenvectors

As with right multiplication of a vector (column vector) with a matrix, we may also consider the left multiplication of a row vector with a matrix. In this case the vector being a row vector is important for the product to be defined. We have the following in this case.

**Definition 9.1** (Left-Eigenvector)**.** Let $A \in \mathbb{C}^{n^2}$ be an $n \times n$ matrix and $\lambda$ be an eigenvalue of $A$. Let $\nu \neq 0$ be an row vector in $\mathbb{R}^n$. Then, $\nu \neq 0$ is called a *left eigenvector* of $A$ associated with the eigenvalue $\lambda$ if $\nu\,(A - \lambda I)\,\nu = 0$.

The right eigenvectors are most common in applications and we normally call them eigenvectors without explicitly referring to right.

We have equivalent counterparts for left-eigenvectors for all the results we mentioned in the previous sections. We will list some of the very important results here.

**Lemma 9.2** (2-Spectrum Lemma for Left-eigenvectors)**.** *If an $n \times n$ diagonalizable matrix $A$ has a spectrum of 2, with $\lambda_1$ and $\lambda_2$ being the two eigenvalues. Then,*

$$(A - \lambda_1 I) \cdot (A - \lambda_2 I) = 0, \tag{9.1}$$

*In addition, the nonzero rows of $A - \lambda_2 I$ are the left-eigenvectors of $A$ for the eigenvalue $\lambda_1$ and vice verse.*

The proof of this result is similar to that of Lemma 3.5, except that we use nonzero rows this time. We can extend the Eigenmatrix Theorem, that is Theorem 3.9 for left-eigenvectors as well and is give below.

**Theorem 9.3** (Eigenmatrix Theorem for Left-eigenvectors)**.** *Let $A$ be a diagonalizable $n \times n$ matrix. Let $\lambda_i, 1 \leq i \leq p \leq n$ be the distinct eigenvalues of $A$. Then,*

$$(A - \lambda_1 I)\,(A - \lambda_2 I) \cdots (A - \lambda_p I) = 0 \tag{9.2}$$

*or*

$$\kappa_{\lambda_1}(A) \cdot \kappa_{\lambda_2}(A) \cdots \kappa_{\lambda_p}(A) = 0 \tag{9.3}$$

*In addition, the nonzero row vectors of the matrix*

$$\kappa_{\lambda_1}(A)\cdots\kappa_{\lambda_{i-1}}(A)\cdot\kappa_{\lambda_{i+1}}(A)\cdots\kappa_{\lambda_p}(A) \tag{9.4}$$

*are the left-eigenvectors of $A$ for the eigenvalue $\lambda_i$ for $i = 2, \ldots, p-1$. The cases for $i = 1$ and $p$ can be considered similarly.*

The proof is similar to that of the Theorem 3.9, except that we consider row vectors in this case.

We do not provide examples for this case since the argument is similar to the right eigenvectors except that in the eigenmatrix theorem, we select rows of the left-most matrix rather than the columns of the right-most matrix. Otherwise the discussion is similar to that of right-eigenvectors.

## 10. Numerical Approach to Find Eigenvectors

Base on the idea of $\kappa$-matrices, we may formulate a new problem to approximate the eigenvector using a minimization approach. The problem can be stated as

**Problem 10.1.** Let $\lambda \neq \mu$ be two eigenvalue of $A$. An approximation to the eigenvector of $A$ associated with $\lambda$ is given by $\nu = x_1a_1 + \cdots + x_pa_p$, where $x = (x_1, \ldots, x_p)$ satisfies the minimization problem

$$\min_{x\in\mathbb{R}^p-\{\mathbf{0}\}} A\left(x_1a_1 + \cdots + x_pa_p\right),$$

in which $a_i$ are the linearly independent column vectors of $\kappa_\mu(A)$.

This is because, according to Corollary 3.11, the eigenvector always lies in the column spaces of the complimentary eigenmatrices. Since 0 is an eigenvalue of $\kappa_\mu(A)$, the columns of $\kappa_\mu(A)$ are linearly independent. Thus, the number of vectors in the linear combination is less than $n$ meaning that the computational time is less than the classical problem that minimizes $(A-\lambda I)x,\ x \in \mathbb{R}^n - \{\mathbf{0}\}$. Numerical evidence shows that that is indeed the case. We leave this problem for future studies.

**Problem 10.2.** There are several computer implementations of the Gauss-Jordan elimination of the eigenmatrices $(A - \lambda I)$ to find the eigenvectors by choosing a basis for the nullspace of $\lambda$. The present method does not require the row-reduction step. Thus, it is an interesting and valuable question to do a complexity analysis of the algorithm of 2-spectrum lemma and the eigenmatrix theorem. Such studies will be done in a future project.

To complete the exposition, we will give a short account to some major applications of eigenvectors.

## 11. Major Applications of Eigenvectors

We only discuss two examples in detail and leave the rest for the discussion section.

**Application 11.1.** As the first application, consider the problem of finding the powers of $A$ for $n \in \mathbb{N}$. If $A$ has a diagonalization given by

$$A = PDP^{-1}$$

where $D$ consists of the eigenvalues of $A$, not necessarily distinct, and $P$ consists of the associated eigenvectors. Then for $n \in \mathbb{N}$, we have

$$A^n = \left(PDP^{-1}\right)\left(PDP^{-1}\right)\cdots\left(PDP^{-1}\right) = PD^nP^{-1}$$

However, finding $D^n$ is easy and straightforward. For example, let

$$D = \begin{pmatrix} \alpha & & & \\ & \beta & & \\ & & \ddots & \\ & & & \gamma \end{pmatrix}$$

Then

$$D^n = \begin{pmatrix} \alpha^n & & & \\ & \beta^n & & \\ & & \ddots & \\ & & & \gamma^n \end{pmatrix}$$

for $n \in \mathbb{N}$. If $A$ is not diagonalizable, we use the Jordan forms to produce similar results.

As the next major application of eigenvectors, consider the problem of *solving a system of differential equations* given below.

**Application 11.2.** First consider the system of first order differential equations:

$$\frac{dx_1}{dt} = a_{11}(t)x_1 + a_{12}(t)x_2 + \cdots + a_{1n}(t)x_n,$$
$$\vdots \qquad\qquad \vdots$$
$$\frac{dx_n}{dt} = a_{n1}(t)x_1 + a_{n2}(t)x_2 + \cdots + a_{nn}(t)x_n.$$

This can be written in the form

$$\mathbf{X} = \begin{pmatrix} x_1(t) \\ \vdots \\ x_n(t) \end{pmatrix}, \quad \text{and} \quad \mathbf{A}(t) = \begin{pmatrix} a_{11}(t) & a_{12}(t) & \cdots & a_{1n}(t) \\ \vdots & \vdots & \ddots & \vdots \\ a_{n1}(t) & a_{n2}(t) & \cdots & a_{nn}(t) \end{pmatrix}.$$

or simply

$$\mathbf{X}' = \mathbf{A}\mathbf{X} \tag{11.1}$$

**Theorem 11.3.** *Let $\lambda_1, \lambda_2, \ldots, \lambda_n$ be $n$ distinct real eigenvalues of the constant coefficient matrix $\mathbf{A}$ of the homogeneous system* (11.1) *and let $\mathbf{V}_1, \mathbf{V}_2, \ldots, \mathbf{V}_n$ be the corresponding eigenvectors. Then the general solution of the homogeneous system on the interval* $(-\infty, \infty)$ *is given by*

$$\mathbf{X} = c_1\mathbf{V}_1 e^{\lambda_1 t} + c_2\mathbf{V}_2 e^{\lambda_2 t} + \cdots + c_n\mathbf{V}_n e^{\lambda_n t}.$$

The cases when the eigenvalues are repeated is complicated and involve generalized eigenvectors. The relevant theory can be found on a standard differential equation book such as [58, p. 345]. We only provide summaries of the existing results here.

11.1. **Eigenvalues with multiplicity more than 1 [58, p.345].** We have the following two cases when an eigenvalue repeats.

Let $A$ be an $n \times n$ matrix and the eigenvalue $\lambda$ has the algebraic multiplicity $m$.

(i) Let an eigenvalue $\lambda$ has the *geometric multiplicity* $m \leq n$. Then there are $m$ linearly independent eigenvectors, $\mathbf{V}_1, \mathbf{V}_2, \ldots, \mathbf{V}_m$ of $A$ corresponding to the eigenvalue $\lambda$. The general solution to the homogeneous system (11.1) is given by

$$\mathbf{X} = c_1\mathbf{V}_1 e^{\lambda t} + c_2\mathbf{V}_2 e^{\lambda t} + \cdots + c_m\mathbf{V}_m e^{\lambda t}.$$

(ii) If there is only one eigenvector (rank 1) associated with the eigenvalue $\lambda$, then there are $m$ *generalized vectors* associated with $\lambda$. In this case, the homogeneous system (11.1) has $m$ linearly independent solution of the form:

$$\begin{aligned} \mathbf{X}_1 &= \mathbf{V}_1 e^{\lambda t} \\ \mathbf{X}_2 &= \mathbf{V}_1 t e^{\lambda t} + \mathbf{V}_2 e^{\lambda t} \\ &\vdots \\ \mathbf{X}_m &= \mathbf{V}_1 \frac{t^{m-1}}{(m-1)!} e^{\lambda t} + \mathbf{V}_2 \frac{t^{m-1}}{(m-1)!} e^{\lambda t} + \cdots + \mathbf{V}_m e^{\lambda t}, \end{aligned}$$

where $\mathbf{V}_i$ are the generalized vectors that satisfy

$$\begin{aligned}(\mathbf{A}-\lambda\mathbf{I})\,\mathbf{V}_1 &= \mathbf{0}\\ (\mathbf{A}-\lambda\mathbf{I})\,\mathbf{V}_2 &= \mathbf{V}_1\\ &\vdots\\ (\mathbf{A}-\lambda\mathbf{I})\,\mathbf{V}_m &= \mathbf{V}_{m-1}.\end{aligned}$$

11.2. **Complex eigenvalues [58, p.350].** We have the following variation when the eigenvalues are complex.

**Theorem 11.4.** *Let $\lambda = \alpha + \beta i$ be a complex eigenvalue of the coefficient matrix $A$ in the homogeneous system* (11.1). *Let $\mathbf{K}$ be the corresponding eigenvector. Then*

$$\begin{aligned}\mathbf{X}_1 &= (\Re(\mathbf{K})\cos\beta t - \Im(\mathbf{K})\sin\beta t)\,e^{\alpha t}\\ \mathbf{X}_2 &= (\Im(\mathbf{K})\cos\beta t + \Re(\mathbf{K})\sin\beta t)\,e^{\alpha t}\end{aligned}$$

*are linearly independent solutions of* (11.1) *on* $(-\infty, \infty)$.

Some other applications of eigenvectors are summarized below. The list is not exhaustive.

11.3. **Further applications of eigenvectors.**

1) **Singular value decomposition for image compression:** The *singular value decomposition* of an $m \times n$ matrix $\mathbf{A}$ is a factorization of the form $\mathbf{U\Sigma V}^*$, where $\mathbf{U}$ is an $m \times m$ *unitary matrix*, $\mathbf{\Sigma}$ is an $m \times n$ matrix with non-negative real numbers on the diagonal, and $\mathbf{V}$ is an $n \times n$ unitary matrix.

   The diagonal entries $\sigma_i = \Sigma_{ii}$ of $\mathbf{\Sigma}$ are known as the *singular values* of $\mathbf{A}$ which are the square roots of the eigenvalues of $AA^*$. The columns of $\mathbf{U}$ and $\mathbf{V}$ are called the *left-singular vectors* and *right-singular vectors* of $\mathbf{A}$, respectively. In addition, the left-singular vectors of $\mathbf{A}$ are a set of *orthonormal eigenvectors* of $AA*$ and the right-singular vectors are a set of orthonormal eigenvectors of $A^*A$.

   In digital image compression techniques, the largest singular values are used to compress an image to a certain compression level. The lesser singular values are used, the higher the compression would be, with the cost of losing the resolution of the image.
2) **Dimensionality Reduction: Principle Component Analysis(PCA):** The principal components correspond to the largest eigenvalues of $A^TA$ and yields the least squared projection onto a smaller dimensional hyperplane in which the *eigenvectors* become the axes of the hyperplane where most of the data is gathered. Dimensionality reduction is used in machine learning and data science as it allows one to understand where most of the variation in the data comes from.
3) **Spectral Clustering:** Clustering is an important aspect of modern data analysis. It allows one to find critical subsystems or patterns within noisy data sets. In particular, *spectral clustering* techniques make use of the spectrum of the similarity matrix of a graph of a network to perform dimensionality reduction before clustering in fewer dimensions. As an important application, the eigenvector of the second smallest eigenvalue of the Laplacian matrix allows one to find the two largest clusters in a network. These are used in social network such as Facebook, Netflix and Youtube to find trends or behavioral patterns for marketing and research purposes.
4) **Low rank factorization for collaborative prediction:** This is used in several sciences to predict certain patters among large networks through an adjacency matrix. It uses the SVD, and consider only the largest eigenvalues of $A^*A$.
5) **Google's PageRank algorithm:** The *PageRank* is an algorithm used in Google search engine to rank webpages according to their importance and it is the first algorithm that was used by the company [39]. It may be considered as one of the historical moments of the importance of eigenvectors in a real-world application. The PageRank algorithm uses the power method to compute the *principal eigenvector* of the adjacency matrix that corresponds to the web-link graph [22,38].

6) **Lorentz transform:** In an extreme situation, light can be considered as an eigenvector of the Lorentz transform [1].
7) **Stability Analysis:** Eigenvalues and eigenvectors are used in stability analysis of systems. See for example [57].
8) **Modal decomposition:** In control theory and dynamical systems, the *modal decomposition* is used to create the dynamic equation for a given system. In mechanical systems such as the one in the introduction, eigenvalues are the natural frequency and the eigenvectors are the mode shapes [51, Section 3.4].

The discussion is inspired in part by [45]. For further discussions about similar applications can be found in the same reference in the StackExchange blog at [45].

## 12. Why Some Matrices are Non-diagonalizable?

The diagonalizability of a matrix depends on whether there are enough vectors to form the transformation matrix $P$. This, in turn, depends on the linear dependency of the column vectors of the corresponding $\kappa$-matrices. In other words, an $n \times n$ matrix is diagonalizable, if it has $n$ linearly independent eigenvectors.

So, the question of why some matrices does not have enough eigenvectors is an interesting question. This topic needs further studies and will be the subject of a future project.

## 13. Adopt in Classroom Instructions

Eigenvectors have several applications in diverse fields such as in engineering, computer science and data science. One of the difficulties we face as faculty when teaching interdisciplinary topics in the classroom is covering the background information and making sure that students are familiar with the prerequisite material and be able to follow the logical flow of the arguments. Most of the times, absence of solid understanding in those material will hinder the knowledge gain that the students make putting the whole teaching-learning process in a jeopardy.

Having a "shortcut approach" to quickly go through the relevant materials may be a positive point, specially when the time constraint is considered. As an experimental approach we will adopt this method for the first time in classroom instructions at Florida Poly University for courses such as Differential Equations, Linear Algebra and Computational Linear Algebra. The outcomes and students surveys will be published in a journal article that focuses on mathematics education.

For example, we also plan to have a collaborative approach to bring the idea with other disciplines as well and share the common pedagogical approach when introducing the concepts. The project involve writing a monograph to cover the preliminaries and the details.

The part of possible courses that may use this approach, outside of the regular linear algebra classes, may include:

— System of differential equations.
— Finding powers of $A$.
— PCA analysis.
— System Analysis.

*Remark* 13.1. "We find the corresponding eigenvectors by finding the nonzero solutions to $(A-\lambda I)\nu = 0$" will be changed to: "We find the eigenvectors by subtracting eigenvalues from the diagonal entries.", if there are only two eigenvalues.

*Remark* 13.2 (Comment about Naming Convention). It is customary to use $\nu$ for eigenvectors, $\lambda$ or $\mu$ for eigenvalues. So the best option for eigenmatrices would be to use $\kappa$. It is also interesting to note that how they are connected to each other. For example, we use matrices ($\kappa$) to generate eigenvalues $(\lambda, \mu)$ and then eigenvalues to find eigenvectors ($\nu$) in that order. They are related from the equation

$$\kappa_\lambda(A)\nu = 0.$$

## 14. List of Shortcut formulas

This is a summary of all the shortcut methods discussed in the present work and will be very useful in applications.

1) For a $2 \times 2$ matrix $A$

**Quick Method of Finding Eigenvectors of a $2 \times 2$ Matrix**

Let $\lambda_1$ and $\lambda_2$ be the distinct eigenvalues of a $2 \times 2$ matrix $A$. The eigenvectors of $A$ can be found by columns of $A$ by subtracting the *complementary eigenvalues* from the diagonal entries, if the resulting columns are nonzero. That is the nonzero columns of the matrix

$$\begin{pmatrix} a-\lambda_2 & b \\ c & d-\lambda_1 \end{pmatrix}$$

are the associated eigenvectors of $\lambda_1$ and $\lambda_2$, respectively.

In the case of repeated eigenvalue, $\lambda$, the nonzero columns of $A - \lambda I$ become the eigenvectors of $A$. If necessary, the resulting columns can be scaled by removing common factors, specially the excessive negative signs, in either case.

2) Common trick for matrices with only two eigenvalues (3.5)

**$2-$ Spectrum Lemma I**

If an $n \times n$ diagonalizable matrix $A$ has a spectrum of 2, with $\lambda_1$ and $\lambda_2$ being two eigenvalues, not necessarily distinct. Then,

$$(A - \lambda_1 I) \cdot (A - \lambda_2 I) = 0, \tag{14.1}$$

In addition, the nonzero columns of $A - \lambda_2 I$ are the eigenvectors of $A$ for the eigenvalue $\lambda_1$ and vice verse. That is, the nonzero columns of the $\kappa$-matrices are the eigenvectors of the complementary eigenvalue.

3) Column Space Approach to Eigenvectors

**Quick Method of Finding Eigenvectors I**

The eigenvectors can be found by a trail-and-error method by guessing the common vector(s) of the column spaces of $\kappa$-matrices of the complementary eigenvalues. For example, to find the eigenvector of $\lambda_1$, find a common vector of $\mathrm{col}(\kappa_{\lambda_2}(E))$ and $\mathrm{col}(\kappa_{\lambda_3}(E))$, if the matrix has three eigenvalues. The method can be easily extended to more than three eigenvalues.

4) Test for Diagonalization of an $n \times n$ Matrix

**Test for Diagonalization of an $n \times n$ Matrix**

Let $A$ be an $n \times n$ matrix $A$. Then $A$ is diagonalizable if and only if the total number of linearly independent vectors $v_i$ produced from the columns of the matrices generated from the products of the form

$$\kappa_{\lambda_1}(A) \cdots \kappa_{\lambda_{i-1}}(A) \cdot \kappa_{\lambda_{i+1}}(A) \cdots \kappa_{\lambda_n}(A) \tag{14.2}$$

for all possible choices of products is $n$.

## Parallel Projects

First, we created a ResearchGate project entitled *Finding Eigenvectors: Fast & Nontraditional Approach* [43] to facilitate a discussion group with both national and international collaborators to continue further development of the subject.

Second, we have also created a supplementary WikiProject for further discussions of the development of the theory. Several illustrative examples are given. Interested readers are referred to [54] and are encouraged to develop the theory there.

The author suggested that Wikipedia may need to come with a new idea called *WikiProjects*, possibly at `https://en.wikipedia.org/WikiProjects` letting the researchers to develop new ideas there, that may lead to breakthrough discoveries one day. `The current Wikipedia policies does not allow the new research findings to appear on there websites`. As an author and a researcher I would like to appear my project at `https://en.wikipedia.org/WikiProjects/Finding_Eigenvectors:_Fast_and_Nontraditional_Method` [55], even though, this is not an option at the time of writing of the paper.

*Remark* 14.1. Among other well-written expository articles, for example, *Math Origins: Eigenvectors and Eigenvalues* by E. R. Tou [50], the excellent collection, *Earliest Known Uses of Some of the Words of Mathematics*, maintained by Jeff Miller [34], the article *Mathematical Words: Origins and Sources* by John Aldrich [3], Steen's article *Highlights in the History of Spectral Theory* and the encyclopedia *MacTutor History of Mathematics Archive* maintained by the School of Mathematics and Statistics at University of St Andrews are excellent sources of information for anyone who is interested in the origin of mathematics.

*Remark* 14.2. Most of the experiments were done on the either Codingground [13] or Symbolab [49]. The former is a great tool to execute MATLAB/Octave Online. Latter is an excellent tool for matrix calculations.

## Conclusion

We have been overlooking the use of columns of a matrix when producing the eigenvectors. The traditional approach involves the Gauss-Jordan elimination of the rows of the characteristic matrix, $(A - \lambda I)$. In the method discussed in this work, we digress from the traditional approach and use only columns to produce the eigenvectors.

We introduced a new concept called the *eigenmatrix* to the already established sequence of terminology: eigenvalues and eigenvectors with the rationale that the eigenmatrix ($\kappa$-matrix) produces eigenvalues ($\lambda$ or $\mu$) and eigenvalues combined with eigenmatrix will produce eigenvectors ($\nu$). With this terminology, we introduce two major results, namely, the *2-Spectrum Lemma* and the *Eigenmatrix Theorem.*

In particular, we show that we may find eigenvectors of a diagonalizable matrix with spectrum 2 just by subtracting the complementary eigenvalue from the diagonal entries if the resulting columns are nonzero. In case the result is a zero column, we select another column to start with.

We also showed that when the spectrum is more than 3, we may produce eigenvectors by a mere matrix-vector multiplication of two $\kappa$-matrices and we referred to this method of finding eigenvectors as the *matrix product approach.*

When the matrix is nondiagonalizable, the vectors produced in this method will correspond to the generalized eigenvectors of the matrix. As a conjecture, we further generalize the Jordan canonical forms for a new class of generalized eigenvectors that are produced by repeated multiples of certain eigenmatrices. We also provide several shortcut formulas to find eigenvectors that does not use echelon forms.

As the final remark, the method can be summarized as "*Finding your puppy at your neighbors'*!" since the eigenvector(s) corresponding to a certain eigenvalue can be found in the column space of any complementary eigenmatrix, an interesting observation.

## Acknowledgement

I would like to specially thank few people during the work of this project, which lasted almost a month. First, my high school mathematics teacher Gunaratne Silva for showing me the beauty of mathematics, one of my colleagues Gayan Wilathgamuwa and my dissertation adviser Jerzy Kocik for their valuable comments and suggestions. I also would like to thank David Knuth for introducing the TEXand to Leslie Lamport for bringing up LATEX. Finally, I would like to thank my wife Pushpika, our son Ryan and our daughter Manisha for allowing me to carry out the project with minimal distraction and lots of love and sacrifices.

Department of Applied Mathematics, Florida Polytechnic University, Lakeland FL 33805, USA.
*E-mail address*: dkatugampola@floridapoly.edu